\documentclass[final]{siamltex}
\usepackage{amsfonts} 
\usepackage{a4wide}
\usepackage{epsfig}
\usepackage{graphicx}
\usepackage{subfig}   
\usepackage[mathcal]{euscript} 
\usepackage{dsfont, pxfonts, verbatim}
\usepackage{mathrsfs}  
\usepackage{stmaryrd}
\usepackage[english]{babel} 
\usepackage{graphicx}
\usepackage[]{caption,subfig}

\newcommand \dps {\displaystyle}
\newcommand \erf {\mbox{erf}} 

\newcommand \be {\begin{equation}}
\newcommand \ee {\end{equation}}
\newcommand \RR		{\mathbb{R}}

\newcommand \del		\partial
\newcommand \eps		\epsilon
\newcommand \calC		{\mathcal{C}}
\newcommand \calH		{\mathcal{H}}
\newcommand \calO		{\mathcal{O}}
\newcommand \calF		{\mathcal{F}}
\newcommand \calP		{\mathcal{P}}
\newcommand \calQ		{\mathcal{Q}}
\newcommand \calL		{\mathcal{L}}
\newcommand \calU		{\mathcal{U}}

\newcommand{\lpcL}		{\gamma_-}
\newcommand{\lpcR}	{\gamma_+}

\newcommand{\fL}		{{f^-}}
\newcommand{\fR}		{{f^+}}

\newcommand{\Id}		{\mathrm{Id}}

\title{Coupling techniques for nonlinear hyperbolic equations. III.
\\
The well--balanced approximation of thick interfaces}
\author{Benjamin Boutin\thanks{Institut de Recherche Math\'ematiques de Rennes, Universit\'e de Rennes I, 263 Avenue du General Leclerc, 35042 Rennes, France. ({\tt benjamin.boutin@univ-rennes1.fr})}, 
\and 
Fr\'ed\'eric Coquel\thanks{Centre de Math\'ematiques Appliqu\'ees \& Centre National de la Recherche Scientifique, Ecole Polytechnique, 91128 Palaiseau, France. ({\tt coquel@cmap.polytechnique.fr)}}, 
\and 
Philippe G. L{\scriptsize e}Floch\thanks{Laboratoire Jacques--Louis Lions \& Centre National de la Recherche Scientifique, Universit\'e Pierre et Marie Curie (Paris 6), 4 Place Jussieu, 75252 Paris, France. {\tt (contact@philippelefloch.org)}} 
\newline 
2000 \textit{AMS Subject Class.} 35L65, 76L05, 76N.
\textit{Keywords and Phrases.} Coupling problem, thick interface model, entropy solution, finite volume approximation, well--balanced, entropy dissipation.
\newline {\tt To appear in:} SIAM Journal of Numerical Analysis (SINUM), 2013. 
}  
 
\begin{document}

\maketitle

\begin{abstract} 
We continue our analysis of the coupling between nonlinear hyperbolic problems across possibly resonant interfaces. In the first two parts of this series, we introduced a new framework for coupling problems which is based on the so--called thin interface model and uses an augmented formulation and an additional unknown for the interface location; this framework has the advantage of avoiding any explicit modeling of the interface structure. In the present paper, we pursue our investigation of the augmented formulation and we introduce a new coupling framework which is now based on the so--called thick interface model. For scalar nonlinear hyperbolic equations in one space variable, we observe that the Cauchy problem is well-posed. Then, our main achievement in the present paper is the design of a new well-balanced finite volume scheme which is adapted to the thick interface model, together with a proof of its convergence toward the unique entropy solution (for a broad class of nonlinear hyperbolic equations). Due to the presence of a possibly resonant interface, the standard technique based on a total variation estimate does not apply, and DiPerna's uniqueness theorem must be used. Following a method proposed by Coquel and LeFloch, our proof relies on discrete entropy inequalities for the coupling problem and an estimate of the discrete entropy dissipation in the proposed scheme.  
\end{abstract}


\section{Introduction}

\subsection{Main objective} 
Mathematical models involving a coupling between distinct nonlinear hyperbolic systems arise in many applications in physics and engineering science and typically involve a non-homogeneous  flux function which exhibits discontinuities in the spatial variable. Over the past decade, a considerable attention has been paid to the so-called `conservative coupling' framework, in which conservation of the unknown is imposed at flux discontinuities \cite{AdimurthiMishraGowda05,AudussePerthame05,BachmannVovelle06,BurgerKarlsen08,BurgerKarlsenTowers09,SeguinVovelle03} (and the references therein). Assuming the physical system to be isolated in the thermodynamical sense, the conservation requirement is nothing but expected, and many problems of interest naturally fall within this framework. In sharp contrast, several other applications of equal importance lead to non-isolated systems, which interplay with (possibly singular) external sources, the latter (on purpose) locally breaking the conservation property. Typical examples are provided by passive or active control devices, while others may fall within this category when understood in a broader sense. Considering, for instance, fluid flow problems, we observe that momentum and/or energy may be locally supplied or tempered by a wide variety of apparatus, ranging from mechanical to electro-magnetical mechanisms. One intends here to minimize singular head loss and pressure drop, or accelerate and heat a gas; these apparatus may also be used for mixing or cooling purposes. We refer for instance to the book by Gad-El-Hak \cite{GadelHak} for a review of current techniques in aerospace and \cite{neptune} for nuclear safety analysis. Mass may be even locally taken from the flow and then injected at a convenient other location \cite{Duret04} in order to prevent oil transportation pipelines from slugging. 

In the design of large systems, the fine scale description of the control is commonly bypassed, and instead the modeling reproduced its net effect as a sharp transition experienced by the flow at the location of the device.  The thermodynamical properties of the flow may be affected by the control, but even if the flux functions are identical, the resulting jump conditions are not in conservation form. Arguments from physics and experiments commonly provide semi-empirical laws which express the right--hand trace at the standing transition as a nonlinear function of the left--hand trace. This function defines the ``transmission conditions'' and provides the basis for a mathematical formulation in terms of  
a kinetic function \cite{LeFloch-ARMA,LeFloch-book} or a family of paths \cite{LeFloch89, DLM}. 

Various {\it ad hoc} numerical methods have been proposed in order to incorporate these transmission conditions. Recent investigations devoted to a finer design of the transient operating conditions of large systems
have revealed several shortcomings  \cite{neptune} and have assessed the need for a rigorous mathematical investigation. 
As far as the problem of coupling fixed interfaces is concerned, the pioneering work \cite{GodlewskiRaviart04} treated the coupling of scalar conservations laws modeling the coupling interface as being infinitely thin. Several extensions to the case of  hyperbolic systems with possibly distinct sizes and involving general transmission conditions have been proposed \cite{AmbrosoEtAl08,BoutinCoquelGodlewski08,BoutinCoquelLeFloch09a}, and transmission conditions are formulated in a weak sense {\it via} a `double' initial value problem (IBVP) formalism. They can be as well understood as a measure source term whose mass precisely defines the expected departure from conservation, as was proposed in \cite{HL} and \cite{LeFloch-ARMA} for nonconservative and interface problems, respectively. Various numerical methods have been proposed in order to exactly capture (isolated) transmission discontinuities in the setting of 
the coupling problem under consideration in the present paper; these methods are called ``well-balanced'' with respect to the singular transmission source term. (See \cite{Amadori08} for related matter and \cite{Coquel12} for a survey.)

A difficulty arising with thin interfaces lies in the fact that the initial value problem, even with apparently well--defined interface conditions, may turn out to be ill--posed, so that the thin interface model does not fully determine the dynamics of the fluid flow. This feature is related to the resonance that may take place at the interface, when waves associated with the fluid have almost vanishing speed. Even in the scalar case \cite{BoutinChalonsRaviart10}, multiple solutions to the initial value problem are exhibited when the coupling interface is resonant. The failure of uniqueness corresponds in fact to a general situation first described in~\cite{IsaacsonTemple92} and further analyzed in~\cite{GoatinLeFloch04}.

The present work is a continuation of our analysis in \cite{BoutinCoquelLeFloch09a,BoutinCoquelLeFloch09b} (to be continued in \cite{BoutinCoquelLeFloch09d}), which is devoted to resonant coupling interfaces. In the first two parts of this series, we introduced a new framework for the mathematical coupling, based on an augmented formulation which has the advantage of avoiding any explicit description of the interface structure and was referred to as the {\sl thin interface model.}  The coupling problem takes
the form of a standard IBVP problem, which can in turn support various regularizing mechanisms. In \cite{BoutinCoquelLeFloch09a}, we relied the self--similar viscosity method by Dafermos \cite{Dafermos73} and established the existence of self-similar solutions (with shock waves) for the coupling problem of two hyperbolic systems (under fairly general assumptions). However, in the limit of vanishing viscosity parameters that we studied in \cite{BoutinCoquelLeFloch09b}, a lack of uniqueness is observed for solutions involving a resonance effect, even in the simple scalar setting. We emphasize that entropy inequalities that would attempt to incorporate at the macroscopic level the fine scale effects modeled by viscous mechanisms, do not lead to a efficient selection principle for thin coupling problems. 

In the present paper, the augmented approach proposed by the authors~\cite{BoutinCoquelLeFloch09a} is shown to lead to another regularization strategy, now based on a {\sl thick interface model,} as we call it.  Roughly speaking, the singular source term modeling the transition is given a smooth profile but the overall regularization technique achieves the key property that the left- and right-hand traces of any isolated transition waves are still exactly captured. Since the source term is entirely localized within the transition wave, it does not act elsewhere and steady solutions of the IBVP problem are thus expected to stay constant outside the transition profile. This assesses the importance of considering isolated transition waves. Importantly, this accuracy property holds for resonant transition waves. It is achieved thanks to a well-balanced strategy that can be traced back to the seminal work by Greenberg and Leroux \cite{GreenbergLeroux92} (see also Bouchut \cite{Bouchut04}, Gosse \cite{Gosse01} and the references therein). This well-balanced property holds for any given regularized profile, in the setting of general pair of fluxes and transmission conditions. The versatility of the method allows us to address a fairly general non-conservative coupling problem for which a given smooth profile may be promoted from experiments and knowledge considerations.

An outline of this paper is as follows. Focusing on equations in one space variable, we briefly recall some of the existing frameworks for the non-conservative coupling of two scalar laws. We then introduce the augmented PDE model with thick coupling interface. Existence and uniqueness of an entropy weak solution follows from the well-known Kru{\v z}kov's  theorem. We propose a scheme for approximating the solution of the Cauchy problem and prove the expected well-balanced property. Due to a reconstruction procedure of the discrete solutions dictated by the well-balanced property, {\it a priori} uniform BV (bounded variation) estimates are not known in general and strong convergence is proved by following the entropy dissipation method of  Coquel and LeFloch~\cite{CoquelLeFloch93}, using DiPerna's uniqueness theorem in the class of entropy measure-valued solutions. Finally, several numerical results illustrate the flexibility of our strategy and its ability to capture various resonant transition waves.

\subsection{Background}
We consider the coupling between two conservation laws at an interface located at the location $x=0$ on the real line $\RR$:
\be
\label{leftrightp}
\del_t w + \del_x f^{\pm}(w) = 0, \quad  \pm x > 0,
\ee
with unknown $w=w(t,x)\in\RR$, defined for $t> 0$ and $x \in \RR$,
where the fluxes $f^{\pm}$ are twice differentiable functions. Prescribing the initial data $w(0,x) = w_0(x)$ at time $t=0$ does not suffice, and an additional condition for modeling the transient exchange of informations at $x=0$ must be supplemented. Such an additional closure is called a {\it coupling condition} and is motivated by the general considerations made in the introduction. Focusing at this stage on piecewise Lipschitz-continuous solutions with bounded left-- and right--hand traces, this coupling condition can be given the form of a nonlinear closure law expressing on the left-hand state as a function of the right-hand state of vice-versa. 
Many different coupling conditions may be introduced and should reflect the precise departure from the conservation property observed in specific applications.
After \cite{AmbrosoEtAl07}--\cite{AmbrosoEtAl07b} and 
\cite{ChalonsRaviartSeguin08}, the coupling condition is formulated from two monotone (increasing, say) functions $\theta_\pm:w\in\RR\mapsto \theta_\pm(w)\in\RR$, so that one writes 
\be
\label{ccstrong}
(\theta_- \circ w)(t,0-) = (\theta_+ \circ w)(t,0+), \qquad t>0.
\ee 
Two approaches are available for general coupling problems, as now discussed. 

On one hand, infinitely thin interface models are based on prescribing coupled boundary conditions, as pioneered by Godlewski and Raviart~\cite{GodlewskiRaviart04}; see
also Boutin et al.~\cite{BoutinChalonsRaviart10}. Boundary conditions are formulated so that (\ref{ccstrong}) is understood in a weak form, following Dubois and LeFloch~\cite{DuboisLeFloch88}.
In order to formulate the Cauchy problem in the half--space $\RR^+$, a boundary condition $b:t\in\RR\mapsto b(t)\in\RR$ is prescribed at $x=0$ and imposed in the sense $w(t,0+)\in\calO_+(b(t))$, where $\calO_+(b(t))= \big\{ \mathcal{W}(0+,b(t),w),w\in\RR \big\}$ denotes the set of admissible traces $\mathcal{W}(0+,b(t),w)$ at $x=0+$ of Riemann solutions 
$\mathcal{W}(\cdot,b(t),w)$
associated with prescribed left--hand state $b(t)$ and arbitrary right--hand state. 
Denoting by $\calO_\pm$ the set of admissible traces at $x=0\pm$ determined from Riemann solutions with flux $f^\pm$ and
 boundary states $\theta_\mp^{-1} \circ \theta_\pm \circ w(t,0\pm)$, we formulate the coupling condition (\ref{ccstrong}) in the weak sense  
$w(t,0\pm) \in\calO_\pm\big( \theta_\pm^{-1} \circ \theta_\pm \circ w(t,0\pm) \big), t>0$. 
It can be checked \cite{BoutinChalonsRaviart10})  
that the proposed weak form 
 for the coupling condition reduces to the strong version (\ref{ccstrong}) as long as a resonance phenomena does not take place at the interface. By resonance, it is meant that waves from either the left-- and/or right--hand problems interact with the interface, so that 
 the continuity property (\ref{ccstrong}) is lost in general and multiple discontinuous solutions may be available. Therefore, an additional selection criterion is required.


\subsection{Augmented model with thin interface}

On the other hand, in the second framework for coupling problems, the thick interface model of the authors \cite{BoutinCoquelLeFloch09a} 
(also considered  earlier in \cite{BoutinCoquelGodlewski08} for the self--similar regularization of scalar equations), we view the coupling interface as a {\sl standing wave for an augmented system} of partial differential equations, by generalizing here the nonconservative reformulation proposed by LeFloch~ \cite{LeFloch88,LeFloch89} for the nozzle flow problem with discontinuous cross section (cf.~also \cite{LT}). 
The standing wave is designed so that  (away from resonance, at least) a complete set of Riemann invariants is available, in agreement with the continuity property (\ref{ccstrong}). 
We propose to consider a new unknown $u=u(t,x)$ defined by ($t>0$) 
\be
\label{u}
u(t,x)= \left\{ 
\begin{array}{ll}
(\theta_- \circ w)(t,x), & x < 0,\\
(\theta_+ \circ w)(t,x), & x>0,
\end{array}
\right.  
\ee
so that the strong coupling condition (\ref{ccstrong}) is equivalent to the continuity condition 
\be
\label{coupling-u}
u(t,0-)=u(t,0+),\quad t>0.
\ee
Let us stress from now on that this convenient reformulation of the coupling condition will play a central role in the derivation of a well-balanced method.
Observe that (\ref{u}) is a well--defined change of variable, since $\del_w\theta_\pm(w)>0$ for all $w\in\RR$. We introduce an unknown $v=v(t,x)$ which coincides with the Heaviside function for all $t>0$:
\be
\label{v}
v(t,x)=v_0(x)=0 \, \quad \mbox{ if } x<0, \quad \qquad
1 \mbox{ if } x>0.
\ee
The value $0$ is meant to recover the given equation in the half line $\RR^-$, while the value 
$1$ represents the given equation in $\RR^+$. Intermediate values of $v$ model a smooth 
transition region from one problem to the other. We thus introduce the augmented model ($t>0$, $x\in\RR$)  and its associated initial condition
\be
 \label{EM} 
\begin{array}{ll}
 \del_t \calC_{0}(u,v) + \del_x \calC_{1}(u,v) - \del_v \calC_{1}(u,v) \del_x v = 0,
\qquad \quad 
 \del_t v = 0,
\\
u_0(x)=\left\{ 
\begin{array}{ll}
(\theta_- \circ w_0)(x),& x<0,
\\
(\theta_+ \circ w_0)(x),& x>0,
\end{array} 
\qquad \qquad \qquad \quad
v_0(x)=\left\{
\begin{array}{ll}
0,& x<0,
\\
1,& x>0. 
\end{array}\right.\right.
\end{array} 
\ee
Here, $w_0$ denotes the initial data for the coupled problem (\ref{leftrightp})
while the functions $\calC_0$ and $\calC_1$ are 
\be
\label{C0} 
\begin{array}{ll}
& \calC_{0}(u,v) = (1-v)\gamma_-(u)+v\gamma_+(u),
\qquad
\quad
\calC_{1}(u,v) = (1-v)f^-(\gamma_-(u))+vf^+(\gamma_+(u)),
\end{array}
\ee 
and $\gamma_\pm$ are defined to be the inverse functions of the increasing map $\theta_\pm$,
respectively. 
By our monotonicity assumption on $\theta_\pm$, one has
\be
\label{monC0}
\del_u\calC_0(u,v)>0,\qquad 
u\in\RR, \quad v\in[0,1].
\ee
This property obviously preserves the time direction determined by  the nonlinear first--order augmented system (\ref{EM}).
Other choices of the coupling functions $\calC_0, \calC_1$ are possible, as we discussed in~\cite{BoutinCoquelLeFloch09a}.
Observe that smooth solutions to (\ref{EM}) obey 
\be
\label{EMnc}
 \del_u \calC_{0}(u,v) \, \del_t u + \del_u \calC_{1}(u,v) \, \del_x u = 0,
\qquad \quad 
 \del_t v = 0, 
\ee
so that the first--order system admits two real eigenvalues: $0$ and 
$\lambda(u,v)=\big( \del_u\calC_0(u,v) \big)^{-1}\del_u\calC_1(u,v)$. 
This system also admits a basis of eigenvectors and the characteristic field associated with 
$\lambda(u,v)$ is genuinely nonlinear, provided the flux functions $f^\pm$ are genuinely nonlinear. The other field is linearly degenerate
and the standing wave is clearly characterized by
$\del_u\calC_1(u,v) \, \del_x u =0$, 
which, provided that $\del_u\calC_1(u,v) \neq 0$, implies the Riemann invariant
$u=\mathrm{cst}$, so that the coupling condition
$u(0-,t)=u(0+,t)$ stated in (\ref{ccstrong}) is satisfied in a strong sense.

States $(u_\star,v_\star)$ with $\del_u\calC_1(u^\star,v^\star)=0$ may exist when one (or both) speed ${(f^\pm)}'$ changes sign. Observe that such states come with the property $\lambda(u^\star,v^\star)=0$ and thus correspond to the interaction of a possibly genuinely nonlinear field with a linearly degenerate one. We must then define weak solutions to the non-conservative nonlinear system (\ref{EMnc}) when eigenvalues vanish. Solutions for such hyperbolic systems are not uniquely defined, unless additional physics is prescribed, as recognized in~LeFloch \cite{LeFloch88,LeFloch89,LeFloch-book}. This situation is referred hereafter as to a resonance phenomena. Indeed it has direct connection with the setting for resonance investigated by Isaacson and Temple~\cite{IsaacsonTemple92} and Goatin and LeFloch~\cite{GoatinLeFloch04}.
As already reported, the definition of weak solutions for (\ref{EM}) in the resonant regime has been tackled in \cite{BoutinCoquelGodlewski08} via the Dafermos self--similar vanishing viscosity analysis.  
In \cite{BoutinCoquelLeFloch09b}, multiplicity of self-similar solutions is shown to persist for the nonconservative model in the limit of vanishing viscosity.
Failure of uniqueness arises for resonance problem as noted in \cite{IsaacsonTemple92} and for various interface problems even in linear hyperbolic equations \cite[Chap.~5]{LeFloch-book}. The origin for multiple solutions 
is found  
in the property that Riemann solutions describe the time--asymptotic behavior of the Cauchy problem for parabolic perturbations of (\ref{EMnc}). Here, the precise definition of a regularization $v_\eta$ (for some $\eta>0$)
plays a central role in the non--uniqueness of Riemann solutions. 
The regularized profile $v_\eta$ does not weight the wave speeds ${f^\pm}'$ equally within the expression
$
\lambda(u,v_\eta)=\big( \del_u\calC_0(u,v_\eta)\big)^{-1} \, 
\big((1-v_\eta)\gamma_-'(u){f^-}'(\gamma_-(u))+v_\eta\gamma_+'(u){f^+}'(\gamma_+(u))\big).
$
Consequently, in the resonance phenomena more importance is given 
to the left or right--hand problem and this is the origin of the failure of uniqueness. We refer the reader to \cite{BoutinCoquelLeFloch09b} where up to four solutions can be build from self--similar analysis. We provide below
 numerical evidences that multiple solutions are stable, in the sense that each solution can be captured numerically. (See also Schecter et al. \cite{LinSchecter03,SchecterPlohrMarchesin04} for a discussion of multiple self-similar solutions.) 


\subsection{Thick interface model} 

We now extend the previous step-like color-function to a {\sl smooth} color-function $v=v(x)$ and we consider 
\be
\label{mapping}
w(t,x)=\calC_0\big(u(t,x),v(x) \big),\qquad t>0,\ x\in\RR.
\ee
In view of the monotonicity property (\ref{monC0}) satisfied by $\calC_0(.,v)$, the function $u$ can be recovered from $w$, the color function $v$ being fixed. With some abuse in the notation, we write 
$w=w(u,v)$ and $u=u(w,v)$. The interest in this change of variable stems from the fact that 
the first equation in (\ref{EMnc}) reduces to the balance law describing the {\sl thick interface model} 
\be
\label{CPw}
\del_t w + \del_x f(w,v) = \ell(w,v) \del_x v, 
\qquad
\quad 
w(0,x)=w_0(x),
\ee
where, in agreement with (\ref{EM}), 
\be
\label{flux}
f(w,v)=\calC_1(u(w,v),v),\qquad \quad \ell(w,v)=\del_v\calC_1(u(w,v),v)
\ee
and the initial data is $w_0=w(u_0,v)$, in agreement with (\ref{EM}). Let us stress that the product $\ell(w,v) \del_x v$ is now nothing but a standard zero-order source term, due to the smoothness of $v$. 

The thick interface framework for coupling problems allows us to use 
the notion of entropy pairs for  conservation laws with source terms. Any convex function $\calU=\calU(w)$ can be used
to define an entropy pair $(\calU,\calF)$. To define the required flux, we start from the augmented system (\ref{EM}) and write, for smooth solutions,
$\del_t \calC_0(u,v(x)) + \del_u\calC_1(u,v(x))\del_xu =0$, which leads us to 
$\del_t \calU(\calC_0(u,v)) + \del_x \calQ(u,v) - \del_v \calQ(u,v)\del_x v=0$, with 
\be
\label{ent1}
\calQ(u,v)=\int^{u}\calU'(\calC_0(\theta,v))\del_u \calC_1(\theta,v)\,d\theta.
\ee
In terms of the unknown $w$, this reads 
$\del_t \calU(w) + \del_x \calF(w,v) - \calL(w,v)\del_xv =0$, 
with
\be
\label{EntFluxw}
\calF(w,v)=\calQ(u(w,v),v),\quad \quad \quad  \calL(w,v)=\del_v\calQ(u(w,v),v).
\ee
Weak solutions of the conservation law with smooth spatial inhomogeneities (\ref{CPw}) are then naturally selected by
the  {\sl entropy inequalities} 
\be
\label{INEw}
\del_t \calU(w) + \del_x \calF(w,v) - \calL(w,v)\del_x v \leq 0, 
\ee
understood in the distributional sense for any convex entropy pair $(\calU,\calF)$. 
Here and since again $v$ is smooth, $\calL(w,v)\del_x v$ acts as a usual source term and
 Kru{\v z}kov's theory \cite{Kruzkov70} applies and provides us with a unique entropy solution to (\ref{CPw})--(\ref{INEw}), when
the flux and source terms are sufficiently regular, say piecewise differentiable. 
The minimal smoothness property on the color function $v$ to meet the Kru{\v z}kov assumptions is therefore 
$v\in W^{2,\infty}(\RR)$. 
We will see that existence and uniqueness of an entropy solution of (\ref{CPw})--(\ref{INEw}) in fact holds 
(under this smoothness condition but) for general initial data, that is, $w_0\in L^\infty(\RR)$.

In order to motivate our method, recall here some properties of {\sl time--independent} solutions to (\ref{CPw}), i.e.~solutions satisfying 
$\del_x f(w,v)=\ell(w,v) \, \del_x v$ 
or, in the $u$--variable, 
$
((1-v)\fL'(\lpcL(u))\lpcL'(u)+v\fR'(\lpcR(u))\lpcR'(u))\del_x u = 0.
$
At the numerical level, it is very challenging to capture steady solutions, especially when the coefficient $((1-v)\fL'(\lpcL)\lpcL'(u)+v\fR'(\lpcR)\lpcR'(u))$ vanishes, that is, when the non--trivial eigenvalue $\lambda(u,v)$ of the hyperbolic system (\ref{EM}) vanishes 
---which is the resonance phenomena. In view of the coupling condition (\ref{coupling-u}), 
our strategy will be to focus on (non constant)
solutions $w$ to (\ref{CPw})--(\ref{INEw}) which 
have {\sl constant component} 
$u$ but {\sl variable component} $v$ which are treated as ``stable solutions'' ---{\sl even when resonance occurs.} 

In the next section, we therefore introduce a {\sl well-balanced finite volume method} for approximating the entropy solution to
 (\ref{CPw})--(\ref{INEw}). As already stressed, by well-balanced we mean that solutions $w(u,v)$ with constant components $u$ and general components $v$, so that $u(w,v)$ is constant in space and time. An adapted {\sl reconstruction procedure} will be required in order to achieve our goal. The resulting family of approximate solutions will be seen to be uniformly bounded in sup--norm under a natural CFL (Courant--Friedrichs--Levy) condition. However a uniform {\it a priori} estimate in total variation is not available (except for the trivial case $w(u,v)=u$), due to the well-balanced reconstruction step.
Following Coquel and LeFloch \cite{CoquelLeFloch93}, we propose to take advantage of the existence of infinitely many entropy differential inequalities and advocate the use of DiPerna's theory of entropy measure-valued solutions \cite{DiPerna85}.


\section{Formulation of the well-balanced scheme and convergence theorem}

\subsection{Notation and assumptions}

In this section, we present a finite volume method for the approximation of entropy solutions satisfying, by definition, (\ref{CPw}) and (\ref{INEw}). As explained, we require the method to be well-balanced with respect to the {\sl family of stationary solutions} $w=w(x)$ relevant for the coupling problem, that is, solutions characterized by the condition that $u(w(x),v(x))$ is constant in $x\in\RR$. Since $v$ depends on $x$, so does the stationary solutions $w(x)$. To ensure the well--balanced property, it is convenient to design a finite volume method and handle the two components of the system on {\sl two distinct grids.} (For this strategy, we refer the reader to the review \cite{Bouchut04}, as well as  \cite{ChainaisCHampier01,EymardGallouetHerbin00,Kroener98}).

For simplicity and without genuine loss of generality, we consider constant time and space steps, denoted by $\Delta t, \Delta x>0$, respectively. We then introduce the time levels $t^n = n\Delta t$ ($n=0,1, \ldots$),
the cell centers $x_j = j\Delta x$, and the cell interfaces $x_{j+1/2} = (j+1/2)\Delta x$ (for all integers $j$). The approximate solutions $u_{\Delta x}$ and $v_{\Delta x}$ are sought as piecewise constant functions, with
\be
\label{delocuv} 
\begin{array}{lll}
u_{\Delta x}(t,x) = u^n_j, \quad &x\in (x_{j-1/2},x_{j-1/2}), \quad &t\in(t^n,t^{n+1}),\\
v_{\Delta x}(t,x) = v_{j+1/2},\quad &x\in(x_j,x_{j+1}), \quad &t\ge 0,
\end{array}
\ee 
and, in view of (\ref{mapping}), we have also the companion function 
\be
\label{defwdx}
w_{\Delta x}(t,x) = w(u_{\Delta x}(t,x), v_{\Delta x}(t,x)),\quad t\ge 0,~x\in \RR.
\ee 
 Since the solution $v$ in (\ref{EM}) is independent of the time variable,  $v_{\Delta x}(t,x)$ is chosen to be time--independent. We also set 
\be
\label{idvdisc}
\begin{array}{lll}
\dps
 u^0_{\Delta x}(x) = u_j^0 = {1 \over \Delta x}\int_{x_{j-1/2}}^{x_{j+1/2}}u_0(x)\,dx, \qquad  x\in (x_{j-1/2},x_{j-1/2}),
\\
\dps
v_{\Delta x}(t,x) 
= v_{j+1/2} = {1 \over \Delta x}\int_{x_j}^{x_{j+1}}v(x)\,dx, \qquad  x\in(x_j,x_{j+1}).
\end{array}
\ee
The discrete solution $u_{\Delta x}$ will be evolved in time by a finite volume method, consistent with the equation of interest 
\be
\label{wlawnum}
\del_t w(u,v) + \del_x f(w(u,v),v) - \big(f_+(\gamma_+(u))-f_-(\gamma_-(u))\big) \, \del_x v = 0.
\ee
By construction, the discrete function $v_{\Delta x}$ is constant within a neighborhood of each cell interface $x_{j+1/2}$, so that the above equation {\sl locally} reduces to a {\sl conservation law} in the unknown $w = w(u,v_{j+1/2})$:
\be
\label{slcloc}
\del_t w + \del_x f(w,v_{j+1/2}) = 0, \quad x\in (x_j,x_{j+1}), \quad t\in (t^n,t^{n+1}).
\ee
This property motivates us to introduce, at each cell interface $x_{j+1/2}$, a two--point numerical flux $g(\cdot,\cdot;v_{j+1/2}) : \RR\times\RR \to \RR$, which is assumed to be locally Lipschitz continuous and satisfy the consistency and monotonicity properties  
\be
\label{consflux}
g(a,a;v_{j+1/2}) =  f(a,v_{j+1/2}),\qquad  a\in \RR
\ee 
\be
\label{monotone}
{\del g \over \del a} (a,b;v_{j+1/2})\ge 0, \qquad 
{\del g \over \del b} (a,b;v_{j+1/2}) \leq 0, 
\qquad 
a,b \in \RR.
\ee 


\subsection{The well-balanced scheme}
\label{sect:wellbalancedscheme}

The discrete solution $u_{\Delta x}(t^n,.)$ being known at time $t^n$, we determine the new approximation at the time $t^{n+1}$ into two steps: a {\sl subcell reconstruction step} followed by an {\sl evolution step.} In turn, our algorithm is a time--explicit finite volume method, which will be shown to converge under the (CFL) stability condition 
\be
\label{cfl}
\frac{\Delta t}{\Delta x} \max_{j} \sup_{u\in [m,M]} \Big|
 \frac{\del f}{\del w}(w(u,v_{j+1/2}),v_{j+1/2})\Big|
 \le \frac{1}{2}, 
\ee
where  $m = \inf_{x\in\RR} u_0(x)$ and $M= \sup_{x\in\RR} u_0(x)$. 
\begin{itemize}
\item{\it Subcell reconstruction.}
At each time $t^n$ and in each cell $(x_{j-1/2},x_{j+1/2})$, using the change of variable $w=w(u,v)$ we introduce two ``subcell states'' $w^n_{j\mp1/2,\pm}$ together with their average: 
\be
\label{subcellw}
\begin{array}{lll}
w^n_{j-1/2,+} = w(u^n_j,v_{j-1/2}),  \quad w^n_{j+1/2,-} = w(u^n_j,v_{j+1/2}),
\qquad
w^n_j = \frac{1}{2}\big(w^n_{j-1/2,+} +w^n_{j+1/2,-} \big).
\end{array} 
\ee
\item{\it Evolution in time.}
At the time $t^{n+1}$ and in each cell $(x_{j-1/2},x_{j+1/2})$, we define $w^{n+1}_j$ by integration on subcells and set  
\be
\label{evolvew}
\begin{array}{lll}
w^{n+1}_j = w^n_j - \frac{\Delta t}{\Delta x}\big(G^n_{j+1/2,-}-G^n_{j-1/2,+}\big), 
\end{array} 
\ee
with $G^n_{j+1/2,-} = g(w^n_{j+1/2,-},w^n_{j+1/2,+};v_{j+1/2}) - f(w^n_{j+1/2,-},v_{j+1/2})$ and 
a similar expression for $G^n_{j-1/2,+}$.  
Then, the ``new state'' $u^{n+1}_j$ is defined as the solution to the algebraic equation 
\be
\label{evolveu}
\frac{1}{2}\big(w(u^{n+1}_j,v_{j-1/2})+w(u^{n+1}_j,v_{j+1/2})\big) = w^{n+1}_j. 
\ee 
\end{itemize}
This completes the description of the proposed method. 

The monotonicity property (\ref{monC0}), namely $\del_uw(u,v)>0$, ensures that
 (\ref{evolveu}) admits a unique solution. Observe that (\ref{evolvew}) is also equivalent to  
$$
\begin{array}{lll}
w^{n+1}_j = w^n_j -\frac{\Delta t}{\Delta x}\Big(g^n_{j+1/2}-g^n_{j-1/2}\Big) 
+\frac{\Delta t}{\Delta x}\big(f(w^n_{j+1/2,-},v_{j+1/2})-f(w^n_{j-1/2,+},v_{j-1/2})\big), 
\end{array}
$$ 
with $g^n_{j+1/2} = g(w^n_{j+1/2,-},w^n_{j+1/2,+};v_{j+1/2})$ and we obtain from the definition (\ref{flux}) of $f(w,v)$
$$
f(w^n_{j+1/2,-},v_{j+1/2})-f(w^n_{j-1/2,+},v_{j-1/2}) = 
 \big(f_+(\gamma_+(u^n_j))-f_-(\gamma_-(u^n_j))\big) (v_{j+1/2}-v_{j-1/2}), 
$$
so that (\ref{evolvew}) is a formally consistent discretization of the governing equation (\ref{wlawnum}) in the unknown $w$.  


\subsection{Main convergence and well--balanced results}

We can now state our main results.

 \vskip.2cm

\begin{theorem}[Convergence of the finite volume method]
\label{thm-convergence}
Consider the Cauchy problem (\ref{EM}) with initial data $u_0$ in $L^\infty(\RR)$ and $v_0$ in $W^{2,\infty}(\RR,[0,1])$
and suppose that the monotonicity property  $\del_u w(u,v)>0$ holds (see (\ref{monC0})). 
Consider the family of functions $v_{\Delta x}$ defined in (\ref{delocuv})--(\ref{idvdisc})
and the family of approximate solutions $w_{\Delta x}$ defined by
 (\ref{defwdx}) (from $u_{\Delta x}$ in (\ref{delocuv})) and the finite volume method (\ref{subcellw})--(\ref{evolvew}) whose
 numerical functions satisfy (\ref{consflux})--(\ref{monotone}). 
Then, under the CFL condition (\ref{cfl}) and as $\Delta x \to 0$, the solutions $w_{\Delta x}$ 
remain bounded in $L^\infty(\RR_+\times\RR)$ and converge strongly in the $L^\infty_t L^p_{loc}$ norm 
($1\le p< \infty$)  toward
 the unique Kru{\v z}kov solution $w$ to (\ref{CPw})--(\ref{INEw}). 
\end{theorem}
 
 \vskip.2cm

\begin{proposition}[Well--balanced property] 
\label{thm-equilibrium} 
Consider the Cauchy problem (\ref{EM}) when the  the initial data $u_0(x) =  u_\star$ ($x\in\RR$) is a constant 
and the data $v_0: \RR \to [0,1]$ is any smooth function. Then, the discrete solution $u_{\Delta x}$ given by
 (\ref{subcellw})--(\ref{evolvew}) is also constant in space, with
$u_{\Delta x}(t^n,x) = u_\star$ ($x\in\RR$) at each time level $t^n$. 
\end{proposition}

 \vskip.2cm

\begin{proof} 
The discrete initial data is now $u_j^0=u_\star$ for all $j$ so, at each interface $x_{j+1/2}$, we get $w_{j+1/2,-}^0=w_{j+1/2,+}^0 = w(u_\star,v_{j+1/2})$, irrespective of the state $v_{j+1/2}$.  
The subscript $\pm$ may be omitted and the flux at $x_{j+1/2}$ reads
$g^0_{j+1/2}=g(w_{j+1/2}^0,w_{j+1/2}^0;v_{j+1/2})=f(w^0_{j+1/2},v_{j+1/2})$ 
(in view of the consistency property (\ref{consflux})). As a consequence, the left-- and right--hand fluxes (\ref{evolvew}) vanish 
identically: $G^n_{j+1/2,-}=G^n_{j+1/2,+}=0$ for all $j$. The scheme (\ref{evolvew}) thus yields $w^{1}_j = w^0_j$ for all $j$ and,
in view of (\ref{subcellw}), $u^1_j$ in each cell satisfies 
$w(u^1_j,v_{j-1/2})+w(u^1_j,v_{j+1/2}) = w(u_\star,v_{j-1/2})+w(u_\star,v_{j+1/2})$. 
In view of the monotonicity property (\ref{monC0}), namely $\del_uw > 0$, we deduce 
$u^1_j = u_\star$, and an induction yields the desired conclusion.
\end{proof}


\subsection{Formulation based on convex combinations}

\begin{figure}[!ht]
\includegraphics[scale=0.85]{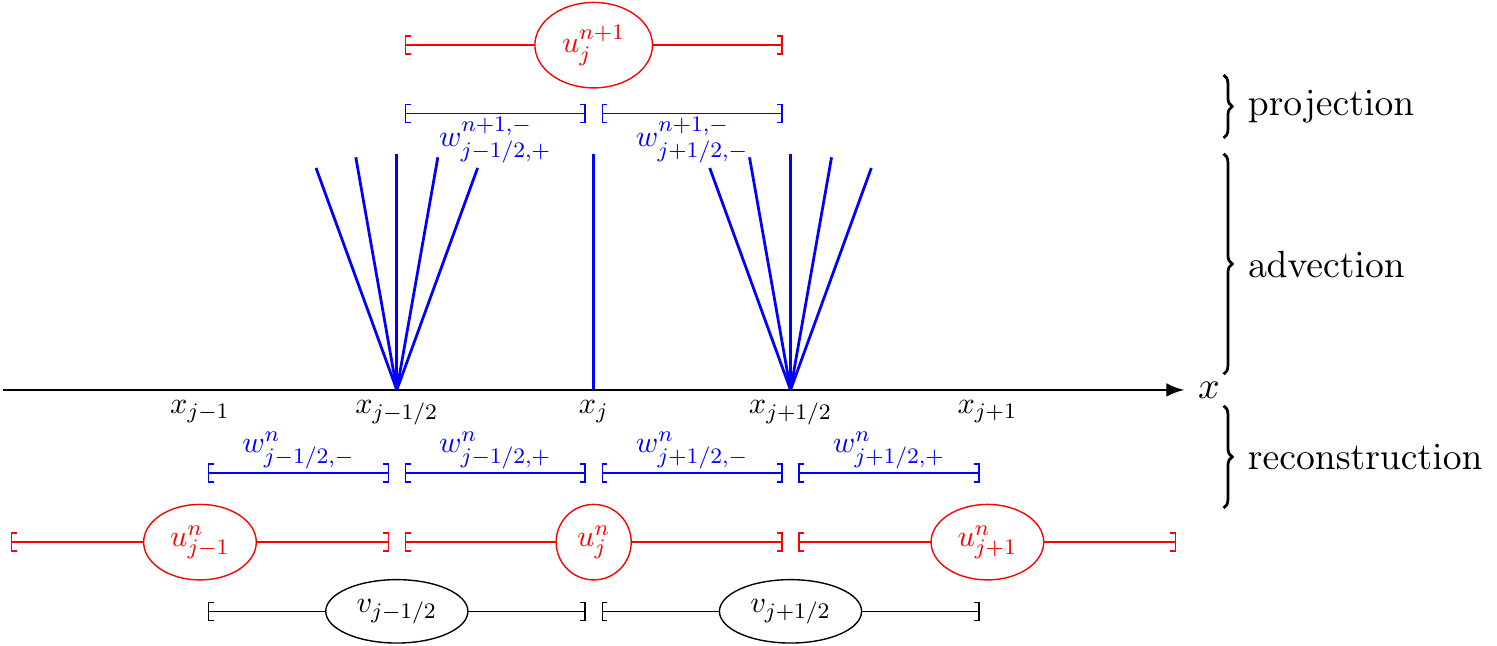}
\caption{A subcell convex combination}
\label{fig-scheme}
\end{figure}

We briefly revisit the finite volume method (\ref{subcellw})--(\ref{evolvew}) so as to highlight its relationships with existing well-balanced approaches. Then, we put forward a convex combination (at the level of subcell) which is of central importance in our forthcoming analysis. 
Consider the following auxilliary Cauchy problem ($t\in(0,\Delta t)$, $x \in \RR$)
\be
\label{couploc}
\del_t w(u,v) + \del_x f(w(u,v),v) - \ell(w(u,v),v)\del_x v = 0,
\qquad
\quad 
\del_t v = 0,
\ee 
\be
\label{riemloc}
\big(u_0(x),v_0(x)\big) = \big(u_{\Delta x}(t^n,x),v^0_{\Delta x}(x)\big) =
\left\{
\begin{array}{lll}
(u^n_j,v_{j-1/2}),\quad x\in (x_{j-1/2},x_j),\\
(u^n_j,v_{j+1/2}),\quad x\in(x_j,x_{j+1/2}). 
\end{array}\right.
\ee
Solving this Cauchy problem in the time slab $(0,\Delta t)$, with $\Delta t$ satisfying the CFL restriction (\ref{cfl}), just amounts to glue together  non--interacting Riemann solutions emanating from the interfaces $x_j$ and $x_{j+1/2}$. Indeed, at $x_j$, the Riemann data has an arbitrary jump in $v_0(x)$ at $x=0$ but with a constant $u_0(x) = u^n_j$, this property allows us to solve the (local) Riemann problem in term of a standing wave
 (depicted  in Figure~\ref{fig-scheme} as a vertical line emanating from $x_j$).  Here, we thus favor such a stationary solution even if resonance locally takes place.  On the other hand, at $x_{j+1/2}$ where $v_0(x)$ is locally constant, the Riemann solution $u$ is 
easily determined at time $\Delta t$, that is, 
$w(u(\Delta t,x),v_{j+1/2})=\omega(\frac{x}{\Delta t};w_{j+1/2,-}^n,w_{j+1/2,+}^n)$ where
$\omega(\frac{x}{t};w_{j+1/2,-}^n,w_{j+1/2,+}^n)$ is the self-similar entropy solution of
\be
\label{lcswloc}
\begin{array}{ll}
&\del_t \omega + \del_x f(\omega,v_{j+1/2}) = 0, \quad t>0, ~x\in\RR, 
\\
&\omega(0,x)=w_0(x) \equiv \left\{
\begin{array}{lll}
w^n_{j+1/2,-} = w(u^n_j,v_{j+1/2}),&x<0,\\
w^n_{j+1/2,+} = w(u^n_{j+1},v_{j+1/2}),&x>0,
\end{array}\right.  
\end{array}
\ee
in agreement with the subcell reconstruction step (\ref{subcellw}). Observe that under the CFL condition (\ref{cfl}), the solution $u(t,x)$ of (\ref{couploc}) cannot interact with the two neighboring standing waves at $x_j$ and $x_{j+1}$ for times $t\in (0,\Delta t)$. Hence the exact solution to the Cauchy problem (\ref{couploc})--(\ref{riemloc}) is obtained in the desired form (Cf.~Fig.~\ref{fig-scheme}).

Next, note that an averaging procedure at time $\Delta t$ in each subcell $(x_{j-1/2},x_j)$ and $(x_j,x_{j+1/2})$ yields
$$
\left.
\begin{array}{lll}
\displaystyle w^{n+1,-}_{j-1/2,+} = \frac{2}{\Delta x}\int_{x_{j-1/2}}^{x_j} w(u(\Delta t,x),v_{j-1/2}) dx, 
\qquad 
\displaystyle w^{n+1,-}_{j+1/2,-} = \frac{2}{\Delta x}\int_{x_j}^{x_{j+1/2}} w(u(\Delta t,x),v_{j+1/2}) dx.
\end{array}
\right.
$$
Using classical arguments, the arithmetic average of the above subcell states gives 
\be
\frac{1}{2}\big( w^{n+1,-}_{j-1/2,+} + w^{n+1,-}_{j+1/2,-}\big) =  w^n_j - \frac{\Delta t}{\Delta x}\Big(G^n_{j+1/2,-}-G^n_{j-1/2,+}\Big) \equiv w^{n+1}_j, 
\ee
where in agreement with (\ref{evolvew})--(\ref{evolveu}), the left-- and right--hand Godunov fluxes read
$G^n_{j+1/2,-} = f(\omega(0+;w^n_{j+1/2,-},w^n_{j+1/2,+}),v_{j+1/2})-f(w^n_{j+1/2,-},v_{j+1/2})$, and
$G^n_{j-1/2,+} = f(\omega(0+;w^n_{j-1/2,-},w^n_{j-1/2,+}),v_{j-1/2})-f(w^n_{j-1/2,+},v_{j-1/2})$.
In other words, the formula (\ref{evolvew}) for $w^{n+1}_j$ is recovered in the special case of the Godunov solver for system (\ref{couploc}). Other left--hand and right--hand fluxes 
(based on general monotone numerical flux functions $g(.,.;v_{j+1/2})$ satisfying (\ref{consflux})--(\ref{monotone}))
are also obtained by approximating the Riemann solution (\ref{lcswloc}).  We summarize our result as follows.

 \vskip.2cm

\begin{lemma}[A subcell convex combination]
\label{lem-scc}
Consider the finite volume method (\ref{subcellw})--(\ref{evolvew}) with fluxes $g(.,.;v_{j+1/2})$ satisfying (\ref{consflux})-(\ref{monotone}), and introduce the subcell states 
\be
\label{eq-updatetime}
 \begin{array}{ll}
  w_{j+1/2,-}^{n+1,-} = w_{j+1/2,-}^n - \frac{2\Delta t}{\Delta x}\left(g_{j+1/2}^n-f(w_{j+1/2,-}^n,v_{j+1/2})\right),\\
  w_{j+1/2,+}^{n+1,-} = w_{j+1/2,+}^n - \frac{2\Delta t}{\Delta x}\left(f(w_{j+1/2,+}^n,v_{j+1/2})-g_{j+1/2}^n\right),
 \end{array}
\ee
with $g^n_{j+1/2}=g(w^n_{j+1/2,-},w^n_{j+1/2,+};v_{j+1/2})$. Then, the formula (\ref{evolvew}) for $w$ recasts in term of the convex combination
\be
\label{convcomb}
w^{n+1}_j = \frac{1}{2}\big(w^{n+1,-}_{j-1/2,+}+w^{n+1,-}_{j+1/2,-}\big).
\ee
\end{lemma}

Rephrasing the above result, the formula (\ref{evolvew}) for $w^{n+1}_j$, which we have underlined to yield a consistent discretization of the equation (\ref{wlawnum}), just reads as a convex combination but for the conservation law (\ref{slcloc}). The next equality reflects such a property and it will be extensively used in the sequel 
\be
\label{idintens}
w^{n+1}_j = \frac{1}{2}\big(w(u^{n+1}_j,v_{j-1/2})+w(u^{n+1}_j,v_{j+1/2})\big) = \frac{1}{2}\big(w^{n+1,-}_{j-1/2,+}+w^{n+1,-}_{j+1/2,-}\big).
\ee


\section{Well-posedness theory for the thick interface model}

To motivate the forthcoming development, we first recall that in the present coupling setting with thick interfaces,  the color function $v$ is a given smooth function of the space variable, say $v\in W^{2,\infty}(\RR)$. Hence as already emphasized, the Kru{\v z}kov's theory applies to establish uniqueness of the entropy solution $w$ in $L^\infty_{loc}(\RR_+\times\RR)$ of the Cauchy problem (\ref{CPw})-(\ref{INEw}) with initial data $w_0\in L^\infty(\RR)$. We will establish hereafter the required sup-norm estimate for the family of approximate solutions $w_{\Delta x}$, with $w_{\Delta x}$ defined in (\ref{defwdx}),   but no such uniform estimate is known in the BV semi-norm. Indeed, the total variation of the discrete solutions may increase at the subcell reconstruction step (\ref{subcellw}) (except in the particular case $\gamma^+=\gamma^-=\Id$, see section~\ref{sec-state}) and a control of the total variation seems out of reach. The absence of an {\it a priori} strong compactness argument leads us to adopt the setting of measured-valued solutions for (\ref{CPw})-(\ref{INEw}) so as to recover, following DiPerna \cite{DiPerna85}, {\it a posteriori} strong convergence of the approximate solutions $w_{\Delta x}$ on the basis of infinitely many entropy inequalities. 

In this section, we state a generalization of DiPerna's uniqueness theorem \cite{DiPerna85} which concerns the class of entropy measure-valued solutions to nonlinear equations (\ref{EMnc}). Measure-valued solutions are Young measures, that is,
 weakly measurable maps $\mu : (t,x) \in\RR_+\times\RR\to \mu_{t,x}$ which take their values in the space of probability measures and, in our case, are supported in a compact interval of $\RR$. A Young measure represent all weak--star limits $a(w_{\Delta x})$ of a bounded sequence $w_{\Delta x}$ for arbitrary functions $a\in\calC^0(\RR)$, that is, 
$
a(w_{\Delta x})
\rightharpoonup
\langle\mu_{t,x}, a\rangle = \int_{\RR_\lambda} a(\lambda) d\mu_{t,x}(\lambda)$
$ \hbox{weakly}-\star \hbox{ in } L^\infty.
$

\begin{definition}[Entropy measure--valued solutions]
\label{def-emvs}
Let $v \in W^{2,\infty}(\RR)$ and $u_0 \in L^\infty(\RR)$ be given, 
and let $w_0=w(u_0,v)$ be an initial data for the Cauchy problem (\ref{CPw}). A measure--valued map $\mu_{t,x}$ is called an entropy measure-valued solution to the Cauchy problem (\ref{CPw}) if for every convex entropy pair $(\calU,{\mathcal F})$ in the form (\ref{ent1})--(\ref{EntFluxw}), one has 
$$
\begin{array}{ll}
\displaystyle \int_{\RR_+\times\RR} \Big(\langle\mu_{t,x}, \calU\rangle\del_t \phi + \langle\mu_{t,x}, {\mathcal F}(\cdot,v)\rangle\del_x \phi \Big) \, dtdx  
 +   \displaystyle \int_{\RR_+\times\RR} \phi \langle\mu_{t,x},{\mathcal L(\cdot,v)}\rangle\del_x v \, dtdx + \int_{\RR} \calU(w_0)\phi(0,.)dx \ge 0 
\end{array} 
$$
for any non--negative test--function $\phi\in{\mathcal D}(\RR_+\times \RR)$.
\end{definition}

With obvious notation, for any continuous function $\calH:\RR\times[0,1]\to\RR$  we set $\langle\mu_{t,x},\calH(\cdot,v)\rangle=\int_\RR\calH(\lambda,v(x))d\mu_{t,x}(\lambda)$. 

 \vskip.2cm

\begin{theorem}[Uniqueness in the class of entropy measure-valued solutions]
\label{theo-uniqemvs}
Let $v$ be given in $W^{2,\infty}(\RR)$, $u_0$ in $L^\infty(\RR)$ and $\mu=\mu_{t,x}$ be an entropy measure-valued solution (in the sense of definition \ref{def-emvs}) of the Cauchy problem (\ref{CPw}) with initial data $w_0=w(u_0,v)$. Then, for almost every $(t,x)$, the measure $\mu_{t,x}$ is a Dirac mass $\mu_{t,x} = \delta_{w(t,x)}$, 
where the function $w\in L^\infty_{loc}(\RR_+\times\RR)$ denotes the unique Kru{\v z}kov's solution of the Cauchy problem (\ref{CPw})--(\ref{INEw}). 
\end{theorem}

 The proof of this result follows from the one in Ben-Artzi and LeFloch in \cite{BenArtziLeFloch07} and 
Amorim, LeFoch, and Okutmustur \cite{ALO}
for conservation laws on manifolds and details are left to the reader.  
Observe here that the initial data is automatically assumed in a strong sense, namely
\be
\label{emvs-id}
\lim_{\tau\to 0\atop \tau>0} \int_0^\tau\int_K \langle\mu_{t,x},|\Id-w_0|\rangle \, dtdx=0
\ee
for all compact $K$ in $\RR$. The following technical lemma provides us with this property, while fturther material can be found in \cite{BenArtziLeFloch07,Szepessy89}. 

\begin{lemma}[DiPerna \cite{DiPerna85}]
\label{lem-idemvs}
Suppose that there exists a strictly convex function $\calU$ and a Young measure $\mu$ satisying, for all $\psi \geq 0$ in  $\calC^\infty_c(\RR)$, 
\be
\label{voirid}
\begin{array}{lll}
\displaystyle 
\lim_{\tau\to 0 \atop \tau>0}\int_0^\tau \int_{\RR} \langle\mu_{t,x}, \Id\rangle \psi \, dtdx = \int_\RR w_0\psi \, dx,
 \qquad 
\lim_{\tau\to 0\atop \tau>0} \int_0^\tau \int_{\RR} \langle\mu_{t,x}, \calU\rangle \psi(x) \, dtdx
 \leq  \int_\RR \calU(w_0)\psi \, dx.  
\end{array}
\ee
Then, the property (\ref{emvs-id}) holds.
\end{lemma}


\section{Convergence analysis}

\subsection{Local maximum principle}

We show now that the approximate solutions remain bounded in $L^\infty(\RR_+\times\RR)$.

 \vskip.2cm

\begin{proposition}[Uniform sup-norm stability]
\label{thm-linfnty}
Under the CFL condition (\ref{cfl}), the finite volume method (\ref{subcellw})--(\ref{evolvew}) satisfies 
\be
\label{ppmu}
\min(u^n_{j-1},u^n_j,u^n_{j+1}) \le u^{n+1}_j\le \max(u^n_{j-1},u^n_j,u^n_{j+1})
\ee
at all time level $t^n$ and, consequently,  
$\| u_{\Delta x}\|_{L^\infty(\RR_+\times\RR)} \le \| u_0\|_{L^\infty(\RR)}$
and $\| w_{\Delta_x}\|_{L^\infty(\RR_+\times\RR)} \le {\mathcal O}(1)$. 
\end{proposition}

{\it Proof.} 
From the subcell states $w^{n+1,-}_{j+1/2,\pm}$ introduced in (\ref{eq-updatetime}), let us define the auxiliary quantity $u^{n+1,-}_{j+1/2,-}$ as the unique solution of
$w^{n+1,-}_{j+1/2,-} \equiv w(u^{n+1,-}_{j+1/2,-},v_{j+1/2})$
as well as $u^{n+1,-}_{j+1/2,+}$ 
from 
$w^{n+1,-}_{j+1/2,+} \equiv w(u^{n+1,-}_{j+1/2,+},v_{j+1/2})$. 
With these definitions, the identity (\ref{idintens}) gives
$w(u^{n+1}_j,v_{j+1/2}) - w(u^{n+1,-}_{j+1/2,-},v_{j+1/2})$
$= w(u^{n+1,-}_{j-1/2,+},v_{j-1/2}) - w(u^{n+1}_j,v_{j-1/2})$. 
We thus deduce 
$
\big(u^{n+1}_j-u^{n+1,-}_{j+1/2,-}\big) \, \big(u^{n+1,-}_{j-1/2,+}-u^{n+1}_j\big) \ge 0 
$
from the monotonicity property satisfied by $w(.,v)$, 
that is
\be
\label{estppmu}
\min(u^{n+1,-}_{j-1/2,+}, u^{n+1,-}_{j+1/2,-}) \le u^{n+1}_j \le \max(u^{n+1,-}_{j-1/2,+}, u^{n+1,-}_{j+1/2,-}).
\ee
It suffices to check that the states $u^{n+1,-}_{j+1/2,\pm}$ satisfy 
$
\min(u^n_j,u^n_{j+1})$ $ \le u^{n+1,-}_{j+1/2,\pm}$ $ \le \min(u^n_j,u^n_{j+1}).
$
For monotone schemes and under the CFL restriction (\ref{cfl}), both subcell states $w^{n+1,-}_{j+1/2,\pm}$ in (\ref{eq-updatetime}) satisfy
\be
\label{ppmw}
\min(w^n_{j+1/2,-},w^n_{j+1/2,+}) \le w_{j+1/2,\pm}^{n+1,-}\le \max(w^n_{j+1/2,-},w^n_{j+1/2,+}). 
\ee
Thanks to the monotonicity of $u(.,v) \equiv w^{-1}(.,v)$ (for $v \in [0,1]$), we deduce that 
$
\min(u^n_j,u^n_{j+1})$ $\le u(w_{j+1/2,\pm}^{n+1,-},v_{j+1/2})$ $\le \max(u^n_j,u^n_{j+1})$, 
since, in the reconstruction step, one has $u^n_j = u(w^n_{j+1/2,-},v_{j+1/2})$ and $u^n_{j+1} = u(w^n_{j+1/2,+},v_{j+1/2})$ in view of (\ref{subcellw}). 
$\hfill\Box$


\subsection{Discrete version of the entropy inequalities}

Since the subcell states $w^{n+1,-}_{j+1/2,\pm}$ are determined from monotone fluxes, they satisfy 
a discrete version of the entropy inequalities coming with (\ref{slcloc}) and we rewrite here $\del_t w + \del_x f(w,v_{j+1/2})=0$. This is the matter of the next statement which is of central importance.

 \vskip.2cm

\begin{lemma}
\label{lem-entnum}
Let $\calU(.), {\mathcal F}(.,v_{j+1/2}) : \RR\to \RR\times \RR$ be any convex entropy pair for the conservation law (\ref{slcloc}). Then there exists a locally Lipschitz continuous, entropy flux ${\mathcal G}(.,.;v_{j+1/2}) : \RR\times\RR\to \RR$
satisfying
 the consistency property ${\mathcal G}(a,a;v_{j+1/2}) = {\mathcal F}(a,v_{j+1/2})$ ($a\in\RR$),  
so that the following discrete entropy inequalities hold  
(under the condition (\ref{cfl}) and with $w_{j+1/2,\pm}^{n+1,-}$ introduced in (\ref{eq-updatetime})):
\be
\label{subcellineq}
\left.
\begin{array}{lll}
\calU(w^{n+1,-}_{j+1/2,-}) -\calU(w^{n}_{j+1/2,-}) + 2\frac{\Delta t}{\Delta x}\big({\mathcal G}^n_{j+1/2}-{\mathcal F}(w^n_{j+1/2,-},v_{j+1/2})\big) \le 0,
\\
\calU(w^{n+1,-}_{j+1/2,+}) -\calU(w^{n}_{j+1/2,+}) + 2\frac{\Delta t}{\Delta x}\big({\mathcal F}(w^n_{j+1/2,+},v_{j+1/2})-{\mathcal G}^n_{j+1/2}\big) \le 0,
\\
{\mathcal G}^n_{j+1/2} = {\mathcal G}(w^n_{j+1/2,-},w^n_{j+1/2,+};v_{j+1/2}). 
\end{array}
\right.
\ee
\end{lemma}

 \vskip.2cm

This result is classical and the proof is omitted. 
The proposed inequalities are the cornerstones of the following (discrete in time but continuous in space) entropy inequalities.

 \vskip.2cm

\begin{lemma}[Shifted time discrete entropy inequalities]
\label{lem-shiftdiscent}
For every test function $\phi \geq 0$ in ${\mathcal D}(\RR^*_+\times\RR)$, one defines  
$\phi^n_{j+1/2} = \frac{1}{\Delta t\Delta x}\int_{t^n}^{t^{n+1}}\!\!\!\int_{x_j}^{x_{j+1}} \phi(t,x)\, dtdx$. 
Then, under the condition (\ref{cfl}), the (discrete--in--time and continuous--in--space) inequality 
\be
\label{contspace}
\left.
\begin{array}{llll}
& \displaystyle \sum_{j} \frac{1}{2}\big(\calU(w^{n+1,-}_{j+1/2,-})+\calU(w^{n+1,-}_{j+1/2,+})\big)\phi^n_{j+1/2}\Delta x
 -\sum_{j} \frac{1}{2}\big(\calU(w^{n}_{j+1/2,-})+\calU(w^{n}_{j+1/2,+})\big)\phi^n_{j+1/2}\Delta x 
\\
& -\displaystyle\int_{t^n}^{t^{n+1}}\!\!\!\int_{\RR} \Big(
 {\mathcal F}(w^n_{\Delta x},v(x))\del_x\phi(t,x) + \phi(t,x) {\mathcal L}(w^n_{\Delta x},v(x))\del_x v\Big) \,  dtdx 
\\
\\
&  \leq
 \displaystyle 
  {\mathcal O}(\Delta x) \Delta t \| \phi \|_{W^{1,\infty}((t^n,t^{n+1})\times \RR)} \chi_{\phi}(t^n,t^{n+1}) 
\end{array}
\right.
\ee
holds, where $\chi_{\phi}(t^n,t^{n+1}) = 1$ if $\max_{ t\in(t^n,t^{n+1})}( \max_{x\in\RR} \phi(t,x)) \not= 0$, while
$\chi_{\phi}(t^n,t^{n+1}) = 0$ otherwise.
\end{lemma}

 \vskip.2cm

In addition to the property $v\in W^{2,\infty}(\RR)$, only the sup--norm estimate in Proposition~\ref{thm-linfnty}   
is needed to deduce  (\ref{contspace}). 
The latter can be handled in the limit $\Delta x \to 0$, by using the Young measure $\mu$ associated with $w_{\Delta x}$. To evaluate the discrete time derivative in (\ref{contspace}), we introduce a cell average representation of $\phi$, i.e. 
\be
\label{cellphi}
\phi^n_j = \frac{1}{2}\big(\phi^n_{j-1/2}+\phi^n_{j+1/2}\big),
\ee 
and we recall 
\be
\label{disctimederiv}
\left.
\begin{array}{lllll}
\displaystyle \sum_j \frac{1}{2}\big(\calU(w^{n+1,-}_{j+1/2,-})\phi^n_{j+1/2}+\calU(w^{n+1,-}_{j-1/2,+})\phi^n_{j-1/2}\big)\Delta x
\displaystyle 
-\sum_j \frac{1}{2}\big(\calU(w^{n}_{j+1/2,-})\phi^n_{j+1/2}+\calU(w^{n}_{j-1/2,+})\phi^n_{j-1/2}\big)\Delta x \\
 \displaystyle = \sum_{j} \Big((\calU(w^{n+1}_j)-\calU(w^{n}_j)\Big) \phi^n_j \Delta x 
  \displaystyle  - \frac{1}{2} \sum_{j} \Big(2~\calU(w^{n+1}_j)\phi^n_j -  \calU(w^{n+1,-}_{j+1/2,-})\phi^n_{j+1/2}-\calU(w^{n+1,-}_{j-1/2,+})\phi^n_{j-1/2}\Big) \Delta x
\\
\quad \displaystyle  + \frac{1}{2} \sum_{j} \Big(2~\calU(w^{n}_j)\phi^n_j -  \calU(w^{n}_{j+1/2,-})\phi^n_{j+1/2}-\calU(w^{n}_{j-1/2,+})\phi^n_{j-1/2}\Big) \Delta x.
\end{array}
\right.
\ee
The first term in the right--hand side of (\ref{disctimederiv}) yields the time derivative, and we control the remaining term as follows.

 \vskip.2cm

\begin{lemma}
\label{lem-errorterm}
The following estimate holds in each mesh cell 
\be
\label{etn} \dps 
2~\calU(w^{n}_j)\phi^n_j -  \calU(w^{n}_{j+1/2,-})\phi^n_{j+1/2}-\calU(w^{n}_{j-1/2,+})\phi^n_{j-1/2}
\le {\mathcal O}((\Delta x)^2) \| \phi \|_{W^{1,\infty}((t^n,t^{n+1})\times(x_{j-1/2},x_{j+1/2}))},
\ee
while
\be
\label{etnpu}
\left.
\begin{array}{lll}
2~\calU(w^{n+1}_j)\phi^n_j -  \calU(w^{n+1,-}_{j+1/2,-})\phi^n_{j+1/2}-\calU(w^{n+1,-}_{j-1/2,+})\phi^n_{j-1/2}
\\
\le -\frac{\sigma_{\calU}}{4} \vert w^{n+1,-}_{j+1/2,-}-  w^{n+1,-}_{j-1/2,+}\vert^2 \phi^n_j 
+ {\mathcal O}(\Delta x) \vert w^{n+1,-}_{j+1/2,-}-  w^{n+1,-}_{j-1/2,+}\vert ~\| \del_x \phi\|_{L^\infty((t^n,t^{n+1})\times(x_{j-1/2},x_{j+1/2}))},
\end{array}
\right.
\ee
where $\sigma_{\calU}$ denotes a convexity modulus of $\calU : \calU''(w) \ge \sigma_{\calU} \geq 0$ for all $w$ such that $|w|\leq\|w_{\Delta x}\|_{L^\infty(\RR^+\times\RR)}$.
\end{lemma}

 \vskip.2cm

With a strictly convex entropy having $\sigma_\calU>0$, the bound (\ref{etnpu}) is slightly sharper than the estimate 
\be
\label{etnpu-bis}
\left.
\begin{array}{lll}
2~\calU(w^{n+1}_j)\phi^n_j -  \calU(w^{n+1,-}_{j+1/2,-})\phi^n_{j+1/2}-\calU(w^{n+1,-}_{j-1/2,+})\phi^n_{j-1/2}
\le  {\mathcal O}(\Delta x) \vert w^{n+1,-}_{j+1/2,-}-  w^{n+1,-}_{j-1/2,+}\vert ~\| \del_x \phi\|_{L^\infty},
\end{array}
\right.
\ee
which is a crucial observation. 
The derivation of (\ref{etnpu}) relies on the following result whose proof is left to the reader.

 \vskip.2cm

\begin{lemma}
\label{lem-refjenest}
Let $w : x\in (0,1) \to w(x)\in\RR$ be  a bounded function with mean value $\overline{w}$. Then, for any convex entropy pair $\calU : w\in \RR \to \calU(w)\in\RR$ the following estimates hold with $m = \min_{x\in (0,1)} w(x)$ and $M = \max_{x\in (0,1)} w(x)$: 
\be
\label{refinedjensenestimate}
\left.
\begin{array}{ll}
\displaystyle \min_{m\le v \le M}\frac{ \calU''(v)}{2}\int^1_0 \vert w(x) -\overline{w}\vert^2 dx 
\le \int_0^1 \calU(w(x))dx - \calU(\overline{w}) 
\leq \dps  \max_{m\le v \le M}\frac{ \calU''(v)}{2}\int^1_0 \vert w(x) -\overline{w}\vert^2 dx.
\end{array}
\right.
\ee
\end{lemma}

Let us comment about (\ref{etn}) and (\ref{etnpu}) in Lemma \ref{lem-errorterm}: with a somehow sharper version of (\ref{etn}) one has 
$$
\left.
\begin{array}{lll}
2~\calU(w^{n}_j)\phi^n_j -  \calU(w^{n}_{j+1/2,-})\phi^n_{j+1/2}-\calU(w^{n}_{j-1/2,+})\phi^n_{j-1/2}
\le \dps
 {\mathcal O}(\Delta x) \vert w^{n}_{j+1/2,-}-  w^{n}_{j-1/2,+}\vert ~\| \del_x \phi\|_{L^\infty((t^n,t^{n+1})\times\RR)}.
\end{array}
\right.
$$
This bound is similar to the estimate (\ref{etnpu-bis}), but has different weights: 
namely,  $\vert w^{n}_{j+1/2,+}-  w^{n}_{j+1/2,-}\vert$ 
instead of $ \vert w^{n+1,-}_{j+1/2,+}-  w^{n+1,-}_{j+1/2,-}\vert$. In the first case, (\ref{subcellw}) for the subcell states $w^{n}_{j+1/2,\pm}$ easily yields the estimate
\be
\label{easyestimate}
\left.
\begin{array}{lll}
 w^{n}_{j+1/2,-}-  w^{n}_{j-1/2,+} 
 =  w(u^{n}_j,v_{j+1/2})- w(u^{n}_j,v_{j-1/2}) 
= {\mathcal O}(1) (v_{j+1/2}-v_{j-1/2}) = {\mathcal O}(\Delta x),
\end{array}
\right.
\ee
as a consequence of Proposition~\ref{thm-linfnty}   for $u_{\Delta x}$ and the smoothness property 
$v\in W^{2,\infty}(\RR)$. By contrast, the jump in the subcell states $\vert w^{n+1,-}_{j+1/2,-} - w^{n+1,-}_{j-1/2,+}\vert$ cannot be expected to vanish uniformly with ${\Delta x}$,
 since discontinuities may develop within the subcells during the evolution step (\ref{eq-updatetime}). One would thus expect to control the jump in (\ref{etnpu}) via a BV estimate, but such an estimate is not available in the present framework. This is the reason why we emphasize the sharper inequality (\ref{etnpu}) over the cruder bound (\ref{etnpu-bis}). The former will be seen to imply the following 
estimate\footnote{Such an estimate was first established, for finite difference schemes in several space dimensions,
 in Coquel and LeFloch~\cite{CoquelLeFloch93}, 
who coined the term ``weak BV estimate'' to denote this bound. The terminology was used in most papers on the subject since then.}.

 \vskip.2cm

\begin{proposition}[Entropy dissipation estimate]  
\label{prop-weakbvestimate}
Let $T>0$ be given and $N_T$ be the greatest integer smaller than $T/\Delta t$. Then, for any time--independent non--negative test function $\psi\in {\mathcal D}(\RR)$, the finite volume approximation (\ref{subcellw})--(\ref{evolvew}) obeys under the CFL condition (\ref{cfl}) the
 weak BV estimate 
\be
\label{weakbvestimate}
\sum_{n=0}^{N_T}\sum_{j} \vert w^{n+1,-}_{j+1/2,-}-w^{n+1,-}_{j-1/2,+}\vert^2 \psi_j \Delta x \ \le {\mathcal O}(1), 
\ee
where  $\psi_j = \frac{1}{2}\big(\psi_{j-1/2}+\psi_{j+1/2}\big)$
and $\psi_{j+1/2} = {1 \over \Delta x}\int_{x_j}^{x_{j+1}} \psi(x)dx$.
\end{proposition}

We will use a slightly simpler estimate, depending upon a non--negative test--function $\psi\in {\mathcal D}(\RR)$ 
such that $\psi(x) = 1,~~\vert x\vert \le L$ for some $L>0$. Denoting $J$ the largest integer smaller than $L/\Delta x$, the estimate (\ref{weakbvestimate}) implies
\be
\label{wbveutile}
\sum_{n=0}^{N_T}\sum_{\vert j\vert < J} \vert w^{n+1,-}_{j+1/2,-}-w^{n+1,-}_{j-1/2,+}\vert^2 \Delta x \ \le {\mathcal O}(1).
\ee


\subsection{Convergence arguments}

The estimate (\ref{wbveutile}) allows us to establish the following (continuous in time and space) version
of the entropy inequalities.

 \vskip.2cm

\begin{proposition}
\label{prop-conttimespace}
Under the CFL condition (\ref{cfl}), the approximate solutions $w_{\Delta x}$ obey the entropy inequality
\[
\left.
\begin{array}{lll}
\displaystyle \int_{\RR_+\times\RR} \Big(\calU(w_{\Delta x})\del_t\phi(t,x) + {\mathcal F}(w_{\Delta x},v)\del_x\phi + \phi {\mathcal L}(w_{\Delta x},v)\del_x v\Big) \, dtdx
+ \int_{\RR} \calU(w^0_{\Delta x})\phi(0,x)dx \ge {\mathcal O}(\sqrt{\Delta x})
\end{array}
\right.
\]
for any (smooth) convex entropy pair $(\calU,\calF)$ in the form (\ref{ent1})--(\ref{EntFluxw}).
\end{proposition}

Equipped with the above inequality,
we deduce that $\mu_{t,x}$ (associated with $\big\{ w_{\Delta x} \big\}_{\Delta x>0}$) satisfies 
\[
\left.
\begin{array}{lll}
\displaystyle \int_{\RR_+\times\RR} \Big(
\langle\mu,\calU(.)\rangle\del_t\phi(t,x) + \langle\mu,{\mathcal F}(.,v)\rangle\del_x\phi + \langle\mu,{\mathcal L}(.,v)\del_x v\rangle\phi \Big) \, dtdx 
+ \int_{\RR} \calU(w_0)\phi(0,x)dx\ge 0.
\end{array}
\right.
\]
In other words, $\mu$ is an entropy measure-valued solution of the Cauchy problem (\ref{CPw}) in the sense of Definition \ref{def-emvs}. Proving that the initial data $\mu_0=\delta_{w_0}$ with $w_0\in L^\infty(\RR)$ the initial data of the problem (\ref{CPw}) is assumed in the strong sense (\ref{voirid}) can be easily deduced from the previous analysis by
following closely related steps (see for instance \cite{CoquelLeFloch93}). The details are left to the reader.
By Theorem \ref{theo-uniqemvs}, the entropy measure--valued solution $\mu_{t,x}$ reduces to
a Dirac measure $\delta_{w(t,x)}$ concentrated on a function $w=w(t, x)$ which coincides with Kru{\v z}kov's solution to
(\ref{CPw})-(\ref{INEw}). This completes the proof of Theorem \ref{thm-convergence}. $\hfill\Box$

\vskip.15cm

{\it Proof of Lemma~\ref{lem-shiftdiscent}.} 
Under the CFL condition (\ref{cfl}), we start from the subcell entropy inequalities (\ref{subcellineq}) satisfied by the subcell states 
$w^{n+1,-}_{j+1/2,-}$ and $w^{n+1,-}_{j+1/2,+}$. Adding these two inequalities yields the following entropy 
inequality centered at $x_{j+1/2}$ 
$$
\begin{array}{lll}
\frac{1}{2}\big(\calU(w^{n+1,-}_{j+1/2,-})+\calU(w^{n+1,-}_{j+1/2,+})\big)
 &- \frac{1}{2}\big(\calU(w^{n}_{j+1/2,-})+\calU(w^{n}_{j+1/2,+})\big) 
\\
&+ \frac{\Delta t}{\Delta x}\big({\mathcal F}(w^{n}_{j+1/2,+},v_{j+1/2})-{\mathcal F}(w^{n}_{j+1/2,-},v_{j+1/2})\big)\le 0. 

\end{array}
$$
Multiplying this inequality by $\phi^n_{j+1/2}\Delta x$ and summing in space gives
\be
\label{discineq}
\begin{array}{lll}
&\sum_{j} \frac{1}{2}\big(\calU(w^{n+1,-}_{j+1/2,-})+\calU(w^{n+1,-}_{j+1/2,+})\big)\phi^n_{j+1/2}\Delta x  
-\sum_{j} \frac{1}{2}\big(\calU(w^{n}_{j+1/2,-})+\calU(w^{n}_{j+1/2,+})\big)\phi^n_{j+1/2}\Delta x \\
 &+ \Delta t \sum_{j} \big({\mathcal F}(w^{n}_{j+1/2,+},v_{j+1/2})-{\mathcal F}(w^{n}_{j+1/2,-},v_{j+1/2})\big)\phi^n_{j+1/2} \le 0. \\
\end{array}
\ee
We now deal with the discrete formulation in space and relying on the identities
$$
\left.
\begin{array}{lll}
\hskip-.22cm  {\mathcal F}(w^{n}_{j+1/2,+},v_{j+1/2}) = 
\frac{1}{2}\big({\mathcal F}(w^n_{j+1/2,+},v_{j+1/2})+{\mathcal F}(w^n_{j+3/2,-},v_{j+3/2})\big)
+\frac{1}{2}\big({\mathcal F}(w^n_{j+1/2,+},v_{j+1/2})-{\mathcal F}(w^n_{j+3/2,-},v_{j+3/2})\big), 
\\ \hskip-.22cm
 {\mathcal F}(w^{n}_{j+1/2,-},v_{j+1/2})= 
\frac{1}{2}\big({\mathcal F}(w^n_{j+1/2,-},v_{j+1/2})+{\mathcal F}(w^n_{j-1/2,+},v_{j-1/2})\big)
+\frac{1}{2}\big({\mathcal F}(w^n_{j+1/2,-},v_{j+1/2})-{\mathcal F}(w^n_{j-1/2,+},v_{j-1/2})\big).
 \end{array}
\right.
$$
We deduce that 
\be
\label{discspaceweak}
\left.
\begin{array}{lll}
\displaystyle \sum_{j}\big({\mathcal F}(w^{n}_{j+1/2,+},v_{j+1/2})-{\mathcal F}(w^{n}_{j+1/2,-},v_{j+1/2})\big)\phi^n_{j+1/2} \\
= 
\displaystyle -\sum_{j}\frac{1}{2}\big({\mathcal F}(w^{n}_{j-1/2,+},v_{j-1/2})+{\mathcal F}(w^{n}_{j+1/2,-},v_{j+1/2})\big)\big(\phi^n_{j+1/2}-\phi^n_{j-1/2}\big) \\
\quad \displaystyle -\sum_{j}  \big({\mathcal F}(w^{n}_{j+1/2,-},v_{j+1/2})-{\mathcal F}(w^{n}_{j-1/2,+},v_{j-1/2})\big)\frac12\big(\phi^n_{j+1/2}+\phi^n_{j-1/2}\big).
\end{array}
\right.
\ee
In view of the definition (\ref{subcellw}) of the subcell states $w^{n}_{j+1/2,\pm}$ and the identity (\ref{EntFluxw}) ${\mathcal F}(w(u,v),v)={\mathcal Q}(u,v)$, we get on one hand 
\be
\label{dswa}
\begin{array}{lll}
&\frac{1}{2}\big({\mathcal F}(w^{n}_{j-1/2,+},v_{j-1/2})+{\mathcal F}(w^{n}_{j+1/2,-},v_{j+1/2})\big) 
 = {\mathcal Q}(u^{n}_j,v_{j-1/2}) \\
&\displaystyle 
+ \frac{1}{2}\int_{0}^{1} \del_v{\mathcal Q}(u^n_j,v_{j-1/2}+s(v_{j+1/2}-v_{j-1/2}))ds\,(v_{j+1/2}-v_{j-1/2})   =  {\mathcal Q}(u^{n}_j,v_{j-1/2}) + {\mathcal O}(\Delta x).
\end{array}
\ee
Here, we have used the sup--norm estimate in Proposition~(\ref{thm-linfnty}), the smoothness of  the mapping ${\mathcal Q}(u,.)$,  together with 
$\vert v_{j+1/2}-v_{j-1/2}\vert = {\mathcal O}(\Delta x)$, 
which is a consequence of the property $v\in W^{2,\infty}(\RR)$
 and (\ref{idvdisc}) defining $v_{\Delta x}(x)$. 

On the other hand, by using similar arguments, we obtain 
\be
\label{dswb}
\left.
\begin{array}{lll}
{\mathcal F}(w^{n}_{j+1/2,-},v_{j+1/2})-{\mathcal F}(w^{n}_{j-1/2,+},v_{j-1/2}) 
\displaystyle 
& \dps
= \int_{0}^{1} \del_v{\mathcal Q}(u^n_j,v_{j-1/2}+s(v_{j+1/2}-v_{j-1/2}))ds\,(v_{j+1/2}-v_{j-1/2}) \\
\displaystyle 
& = \del_v{\mathcal Q}(u^n_j,v_{j-1/2}) (v_{j+1/2}-v_{j-1/2})  + {\mathcal O}((\Delta x)^2).
\end{array}
\right. 
\ee
Plugging (\ref{dswa}) and (\ref{dswb}) in (\ref{discspaceweak}) then gives us 
\be
\label{discspaceweakb}
\left.
\begin{array}{lll}
\displaystyle \sum_{j}\big({\mathcal F}(w^{n}_{j+1/2,+},v_{j+1/2})-{\mathcal F}(w^{n}_{j+1/2,-},v_{j+1/2})\big)\phi^n_{j+1/2} 
\\
= \displaystyle -\sum_{j} {\mathcal Q}(u^n_j,v_{j-1/2}) \frac{\big(\phi^n_{j+1/2}-\phi^n_{j-1/2}\big)}{\Delta x} \Delta x 
\displaystyle -\sum_{j}  \frac12\big(\phi^n_{j+1/2}+\phi^n_{j-1/2}\big) \del_v{\mathcal Q}(u^n_j,v_{j-1/2}) \frac{\big(v_{j+1/2}-v_{j-1/2}\big)}{\Delta x} \Delta x 
\\
\quad \displaystyle + {\mathcal O}({\Delta x}) \| \phi\|_{W^{1,\infty}(\RR_+\times \RR)} \chi_\phi(t^n,t^{n+1}).
\end{array}
\right.
\ee
Finally, routine arguments based on the smoothness of $\phi$ and $v$ yield the expected result 
\be
\label{discspaceweakc}
\left.
\begin{array}{lll}
\displaystyle \Delta t \sum_{j}\big({\mathcal F}(w^{n}_{j+1/2,+},v_{j+1/2})-{\mathcal F}(w^{n}_{j+1/2,-},v_{j+1/2})\big)\phi^n_{j+1/2} 
\\
= \displaystyle -\int_{t^n}^{t^{n+1}}\int_{\RR} \Big(
{\mathcal Q}(u_{\Delta x},v(x)) \del_x\phi \, dtdx +\phi \del_v{\mathcal Q}(u_{\Delta x},v(x)) \del_x v\Big) \, dtdx 
 \\
\quad \displaystyle + {\mathcal O}({\Delta x}) \Delta t \| \phi\|_{L^\infty(\RR_+\times \RR)}\chi_\phi(t^n,t^{n+1}),
\end{array}
\right.
\ee
where, by definition, (\ref{ent1})-(\ref{defwdx}), one has ${\mathcal F}(w_{\Delta x},v(x))={\mathcal Q}(u_{\Delta x},v(x))$ and  ${\mathcal L}(w_{\Delta x},v(x))=\del_v{\mathcal Q}(u_{\Delta x},v(x))$.  

\vskip.15cm 

{\it Proof of Lemma~\ref{lem-errorterm}.} 
We derive (\ref{etn}) by plugging $\big(2~\calU(w^{n}_j)\phi^n_j-\calU(w^{n}_{j+1/2,-})\phi^n_{j+1/2}-\calU(w^{n}_{j-1/2,+})\phi^n_{j-1/2}\big)$ in 
$$
\left.
\begin{array}{lll}
& \calU(w^{n}_{j+1/2,-})\phi^n_{j+1/2}+\calU(w^{n}_{j-1/2,+})\phi^n_{j-1/2} 
\\
&=\big(\calU(w^n_{j+1/2,-})+ \calU(w^n_{j-1/2,+})\big) \phi^n_j 
+ \frac{1}{2}\big(\calU(w^n_{j+1/2,-})-\calU(w^n_{j-1/2,+})\big)\big(\phi^n_{j+1/2}-\phi^n_{j-1/2}\big),
\end{array}
\right.
$$
with $\phi_j^n$ defined in (\ref{cellphi}).
In view of  (\ref{subcellw}) (namely, $2 w^n_j = w^n_{j-1/2,+}+w^n_{j+1/2,-}$),
Jensen's inequality implies 
\be
\label{etnb}
\left.
\begin{array}{lll}
2~\calU(w^{n}_j)\phi^n_j-\calU(w^{n}_{j+1/2,-})\phi^n_{j+1/2}-\calU(w^{n}_{j-1/2,+})\phi^n_{j-1/2} 
\\
= \big(2~\calU(w^{n}_j)-\calU(w^n_{j+1/2,-})- \calU(w^n_{j-1/2,+})\big) \phi^n_j 
+ \frac{1}{2} \big(\calU(w^n_{j-1/2,+})-\calU(w^n_{j+1/2,-})\big)\big(\phi^n_{j+1/2}-\phi^n_{j-1/2}\big) \\
 \le \frac{1}{2}\big(\calU(w^n_{j-1/2,+})-\calU(w^n_{j+1/2,-})\big)\big(\phi^n_{j+1/2}-\phi^n_{j-1/2}\big),
\end{array}
\right.
\ee
which provides the required bound (\ref{etn}) in view of the estimate (\ref{easyestimate}) satisfied by
 the jump $(w^n_{j+1/2,-} - w^n_{j-1/2,+})$. 
Deriving the second estimate (\ref{etnpu}) is completely analogous, from an identity similar to (\ref{etnb}), that is 
\be
\label{etnc}
\left.
\begin{array}{lll}
2~\calU(w^{n+1}_j)\phi^n_j-\calU(w^{n+1,-}_{j+1/2,-})\phi^n_{j+1/2}-\calU(w^{n+1,-}_{j-1/2,+})\phi^n_{j-1/2} 
\\
= \big(2~\calU(w^{n+1}_j)-\calU(w^{n+1,-}_{j+1/2,-})- \calU(w^{n+1,-}_{j-1/2,+})\big) \phi^n_j \\
\quad+ \frac{1}{2} \big(\calU(w^{n+1,-}_{j-1/2,+})-\calU(w^{n+1,-}_{j+1/2,-})\big)\big(\phi^n_{j+1/2}-\phi^n_{j-1/2}\big).
\end{array}
\right.
\ee
In view of the identity $2w^{n+1}_j = w^{n+1,-}_{j+1/2,-}+w^{n+1,-}_{j-1/2,+}$ in  (\ref{convcomb}) Lemma \ref{lem-scc}, we 
can apply Lemma \ref{lem-refjenest} and obtain 
$$
2~\calU(w^{n+1}_j)-\calU(w^{n+1,-}_{j+1/2,-})- \calU(w^{n+1,-}_{j-1/2,+}) \le - \frac{1}{4}\sigma_\calU \vert w^{n+1,-}_{j+1/2,-}-w^{n+1,-}_{j-1/2,+}\vert^2,
$$
so that (\ref{etnpu}) follows from (\ref{etnc}).
This concludes the proof of Lemma~\ref{lem-errorterm}.  $\hfill\Box$

\vskip.15cm

{\it Proof of Proposition \ref{prop-weakbvestimate}.} 
We start from the discrete in time, continuous in space formulation (\ref{contspace}) 
and write 
\be
\label{contspacebis}
\left.
\begin{array}{lll}
& \displaystyle \sum_{j} \frac{1}{2}\big(\calU(w^{n+1,-}_{j+1/2,-})+\calU(w^{n+1,-}_{j+1/2,+})\big)\phi^n_{j+1/2}\Delta x 
 -\sum_{j} \frac{1}{2}\big(\calU(w^{n}_{j+1/2,-})+\calU(w^{n}_{j+1/2,+})\big)\phi^n_{j+1/2}\Delta x 
\\
& -\displaystyle\int_{t^n}^{t^{n+1}}\!\!\!\int_{\RR} \Big(
 {\mathcal F}(w^n_{\Delta x},v(x))\del_x\phi(t,x) + \phi(t,x) {\mathcal L}(w^n_{\Delta x},v(x))\del_x v\Big) \,  dtdx 
\\
 & \le {\mathcal O}(\Delta x) \Delta t \| \phi \|_{W^{1,\infty}((t^n,t^{n+1})\times \RR)} \chi_{\phi}(t^n,t^{n+1}),
\end{array}
\right.
\ee
in which we plug the decomposition (\ref{disctimederiv}) of the discrete time derivative. Here, we use the 
discrete test function $\psi_j$ defined in the proposition
 from any time-independent test function $\psi\in{\mathcal D}(\RR)$. We thus get
\be
\label{contspaceter}
\left.
\begin{array}{lll}
& \displaystyle \sum_{j} \big(\calU(w^{n+1}_j)-\calU(w^{n}_j)\big)\psi_j\Delta x 
  -\displaystyle\int_{t^n}^{t^{n+1}}\!\!\!\int_{\RR} \Big(
 {\mathcal F}(w^n_{\Delta x},v(x))\del_x\psi(x) + \psi(x) {\mathcal L}(w^n_{\Delta x},v(x))\del_x v\Big) \,  dtdx 
\\
&\leq \displaystyle \frac{1}{2} \sum_{j} \Big(2~\calU(w^{n+1}_j)\psi_j -  \calU(w^{n+1,-}_{j+1/2,-})\psi_{j+1/2}-\calU(w^{n+1,-}_{j-1/2,+})\psi_{j-1/2}\Big) \, \Delta x
\\
& \displaystyle  - \frac{1}{2} \sum_{j} \Big(2~\calU(w^{n}_j)\psi_j -  \calU(w^{n}_{j+1/2,-})\psi_{j+1/2}-\calU(w^{n}_{j-1/2,+})\psi_{j-1/2}\Big) \,  \Delta x 
 + {\mathcal O}(\Delta x) \Delta t \| \psi \|_{W^{1,\infty}(\RR)}
\end{array}
\right.
\ee
and the estimates (\ref{etn}) and (\ref{etnpu}) stated in Lemma \ref{lem-errorterm} then yield
\be
\label{contspacequadro}
\left.
\begin{array}{lll}
& \displaystyle \sum_{j} \big(\calU(w^{n+1}_j)-\calU(w^{n}_j)\big)\psi_j\Delta x 
+\frac{1}{8}\sigma_\calU\sum_{j} \vert w^{n+1,-}_{j+1/2,-}-w^{n+1,-}_{j-1/2,+}\vert^2 \psi_j \Delta x 
\\
&\leq \displaystyle  {\mathcal O}(\Delta x) \Delta t \| \psi \|_{W^{1,\infty}(\RR)}
 + {\mathcal O}(\Delta x) \sum_{j} \|\del_x\psi\|_{L^\infty((x_{j-1/2},x_{j+1/2}))} \Delta x
\displaystyle
 + {\mathcal O}((\Delta x)^2) \sum_{j} \|\del_x\psi\|_{W^{1,\infty}((x_{j-1/2},x_{j+1/2}))} \Delta x
\\
&\quad
    +\displaystyle\int_{t^n}^{t^{n+1}}\int_{\RR} \Big(
 {\mathcal F}(w^n_{\Delta x},v(x))\del_x\psi(x) + \psi(x) {\mathcal L}(w^n_{\Delta x},v(x))\del_x v\Big)  dtdx.
\end{array}
\right.
\ee
The sup-norm estimate in Proposition~\ref{thm-linfnty}  implies the following crude estimate
\[
\displaystyle\int_{t^n}^{t^{n+1}}\int_{\RR} \Big(
 {\mathcal F}(w^n_{\Delta x},v(x))\del_x\psi(x) + \psi(x) {\mathcal L}(w^n_{\Delta x},v(x))\del_x v\Big) \,  dtdx 
= {\mathcal O}(\Delta t) \|\psi\|_{W^{1,\infty}(\RR)}, 
\]
so that  
\be
\label{contspacequinte}
\left.
\begin{array}{lll}
\displaystyle \sum_{j} \big(\calU(w^{n+1}_j)-\calU(w^{n}_j)\big)\psi_j\Delta x 
+\frac{1}{8}\sigma_\calU\sum_{j} \vert w^{n+1,-}_{j+1/2,-}-w^{n+1,-}_{j-1/2,+}\vert^2 \psi_j \Delta x 
\le {\mathcal O}(\Delta t) \| \psi \|_{W^{1,\infty}(\RR)}.
\end{array}
\right.
\ee
Summing over all $n\in[0,N_T]$ with $N_T=[T/\Delta t]$ (for $T>0$ fixed), we get
$$
\int_{\RR}\calU(w_{\Delta x}(x,T))dx + \frac{1}{8}\sigma_\calU\sum_{n=0}^{N_T}\sum_{j} \vert w^{n+1,-}_{j+1/2,-}-w^{n+1,-}_{j-1/2,+}\vert^2 \psi_j \Delta x 
\le \displaystyle \int_{\RR}\calU(w_0(x))dx + {\mathcal O}(1) T \| \psi \|_{W^{1,\infty}(\RR)},
$$
which gives the desired estimate (\ref{weakbvestimate}), choosing for instance the quadratic entropy $\calU(w)=w^2/2$ with $\sigma_\calU=1$.
This completes the proof of Proposition \ref{prop-weakbvestimate}. $\hfill\Box$

\vskip.15cm 

{\it Proof of Proposition \ref{prop-conttimespace}.}
For any $\phi\in{\mathcal D}(\RR_+\times\RR)$ and in view of its discrete representation $\phi^n_j$, we again consider the continuous in space formulation (\ref{contspace}) in Lemma \ref{lem-shiftdiscent}. Plugging in the decomposition (\ref{disctimederiv}) gives 
\[
\left.
\begin{array}{llll}
& \displaystyle \sum_{j} \big(\calU(w^{n+1}_j)-\calU(w^{n}_j)\big)\phi^n_j\Delta x  
  -\displaystyle\int_{t^n}^{t^{n+1}}\int_{\RR} \Big(
 {\mathcal F}(w^n_{\Delta x},v(x))\del_x\phi(t,x) + \phi(t,x) {\mathcal L}(w^n_{\Delta x},v(x))\del_x v\Big) \,  dtdx 
\\
&\leq \displaystyle  {\mathcal O}(\Delta x) \Delta t \| \phi \|_{W^{1,\infty}((t^n,t^{n+1}\times\RR)}\chi_\phi(t^n,t^{n+1})
  + {\mathcal O}((\Delta x)^2) \sum_{j} \|\del_x\phi\|_{W^{1,\infty}((t^n,t^{n+1})\times(x_{j-1/2},x_{j+1/2}))} \Delta x
\\
&\quad     \displaystyle + {\mathcal O}(\Delta x) \sum_{j} \vert w^{n+1,-}_{j+1/2,-}-w^{n+1,-}_{j-1/2,+}\vert ~\|\del_x\phi\|_{L^\infty((t^n,t^{n+1})\times(x_{j-1/2},x_{j+1/2}))} \Delta x,
 \end{array}
\right.
\]
where we have used the estimates (\ref{etn}) and (\ref{etnpu-bis}). Summing this inequality over the time indices yields
\[
\left.
\begin{array}{llll}
&\displaystyle -\sum_{n=0,1,\ldots}\sum_{j} \calU(w^{n+1}_j)\frac{\phi^{n+1}_j-\phi^n_j}{\Delta t}\Delta t \Delta x  
  -\displaystyle\int_{\RR_+}\int_{\RR} \Big(
{\mathcal F}(w^n_{\Delta x},v(x))\del_x\phi(t,x) + \phi(t,x) {\mathcal L}(w^n_{\Delta x},v(x))\del_x v\Big) \,  dtdx 
\\
& \leq 
\displaystyle  {\mathcal O}(\Delta x) \| \phi \|_{W^{1,\infty}(\RR_+\times\RR)}
\\
& \quad + \displaystyle {\mathcal O}(1) \sum_{n\ge 0}\sum_{j} \Big(\vert w^{n+1,-}_{j+1/2,-}-w^{n+1,-}_{j-1/2,+}\vert \chi_{\phi}(t^n,t^{n+1})~\|\del_x\phi\|_{L^\infty((t^n,t^{n+1})\times(x_{j-1/2},x_{j+1/2}))} \Delta t\Delta x\Big).
 \end{array}
\right.
\]
Cauchy-Schwarz inequality then allows us to upper bound the last term according to
\[
\left.
\begin{array}{lll}
& \displaystyle\sum_{n\ge 0}\sum_{j} \Big(\vert w^{n+1,-}_{j+1/2,-}-w^{n+1,-}_{j-1/2,+}\vert \chi_{\phi}(t^n,t^{n+1})\Big)~\|\del_x\phi\|_{L^\infty} \Delta t\Delta x 
\\
& \leq
\displaystyle \Big(\sum_{n\ge 0}\sum_{j} \big(\vert w^{n+1,-}_{j+1/2,-}-w^{n+1,-}_{j-1/2,+}\vert \chi_{\phi}(t^n,t^{n+1})\big)^2\Delta t\Delta x\Big)^{1/2}  
 \Big(\sum_{n\ge 0}\sum_{j} \big(\|\del_x\phi\|_{L^{\infty}((t^n,t^{n+1})\times (x_{j-1/2},x_{j+1/2})}\big)^2\Delta t\Delta x\Big)^{1/2} 
\\
& \leq
\displaystyle {\mathcal O}(1)  \, \Big(\sum_{n\ge 0}\sum_{j} \vert w^{n+1,-}_{j+1/2,-}-w^{n+1,-}_{j-1/2,+}\vert^2 \chi_{\phi}(t^n,t^{n+1})\Delta t\Delta x\Big)^{1/2}. 
\end{array}
\right.
\]
The weak BV estimate (\ref{wbveutile}) implies  
that $\sum_{n\ge 0}\sum_{j} \vert w^{n+1,-}_{j+1/2,-}-w^{n+1,-}_{j-1/2,+}\vert^2 \chi_{\phi}(t^n,t^{n+1})\Delta t\Delta x$
is of order ${\mathcal O}(\Delta x)$, at most, 
and routine arguments allow us to conclude. $\hfill\Box$


\section{Coupling via conservative variable and total variation estimate}
\label{sec-state}

This section briefly addresses the particular case of the state coupling, namely the case of the coupling condition (\ref{ccstrong}) with $\theta_-=\theta_+=\Id$, that is, we impose here 
$w(0_-,t)=w(0_+,t)$ for all $t>0$, as considered by Godlewski and Raviart \cite{GodlewskiRaviart04}. 
 In this setting, we can derive a BV estimate for the approximate solutions $w_{\Delta x}$ and the convergence follows from Helly's compactness theorem. No smoothness is required on the color function $v$, 
which can be {\sl discontinuous,} so that the result in this section holds in a stronger norm 
than the one in Theorem~\ref{thm-convergence}, but applies to initial data $w_0\in L^\infty(\RR)\cap BV(\RR)$, only. 
To avoid technicalities, we restrict attention to the Godunov solver. Note that for the coupling  $w(0_-,t)=w(0_+,t)$, one has
$u(w,v)=w$ for all $v$, and the condition (\ref{sc-cfl}) below 
yields $|w|\leq\|w_0\|_{L^\infty(\RR)}$: in view of Proposition~\ref{thm-linfnty}, one has $\|w\|_{L^\infty(\RR^*_+\times\RR)}\leq \|w_0\|_{L^\infty(\RR)}$, while $\del_wf(w,v)$ coincides with the non-trivial eigenvalue of (\ref{EM}).

\vskip.20cm 

\begin{theorem}[Total variation estimate]
\label{thm-BV}
Assume $\theta_-=\theta_+=\Id$. Let $v$ in $L^\infty(\RR)$ and $w_0$ in $L^\infty(\RR)\cap BV(\RR)$ be the initial data. For any small $\epsilon>0$ ($\epsilon<1/2$) and under the strengthened CFL condition
\be
\label{sc-cfl}
\max_{|w|\leq \|w_0\|_{L^\infty(\RR)}\atop |v|\leq 1} |\del_wf(w,v)| \leq \frac12-\epsilon,
\ee
the finite volume method (\ref{subcellw})--(\ref{evolvew}) is total variation diminishing and, in particular, 
$TV(w_{\Delta x}(t^{n},\cdot))\leq TV(w_0)$. 
\end{theorem}

{\it Proof.} 
The case of the state coupling comes with $\gamma_-=\gamma_+=\Id$, therefore the reconstruction step boils
 down to $w_{j-1/2,+}^n=w_{j+1/2,-}^n=w_j^n$, namely the discrete solution is kept unchanged and no jump is created at this step. 
Consider the Riemann problem at the interface $x_{j+1/2}$ (with $v_{j+1/2}$ constant) and denote $w(\cdot,w_j^n,w_{j+1}^n)$ its self-similar solution, which is known to be monotone and satisfies 
$TV(w(.,w_j^n,w_{j+1}^n)) = |w_{j+1}^n-w_j^n|$. 
Define $w_{\Delta x}(t^{n+1,-},x)$ by glueing together non-interacting neighboring Riemann solutions. Under the CFL condition (\ref{sc-cfl}), we find
$$
\begin{array}{lll}
 TV(w_{\Delta x}(t^{n+1-},\cdot))
 & = \sum_j TV_{[x_j+\epsilon\Delta x/2,x_{j+1}-\epsilon\Delta x/2]}(w_{\Delta x}(t^{n+1-},\cdot))
 + \sum _j TV_{(x_j-\epsilon\Delta x/2,x_j+\epsilon\Delta x/2)}(w_{\Delta x}(t^{n+1-},\cdot))\\
 & = \sum_j TV(w(\cdot,u_j^n,w_{j+1}^n)) + \sum _j TV_{(x_j-\epsilon\Delta x/2,x_j+\epsilon\Delta x/2)}(w_{\Delta x}(t^{n},\cdot)),
\end{array}
$$
but the second term vanishes since $w_{j-1/2+}^n=w_{j+1/2-}^n$. Therefore, 
$TV(w_{\Delta x}(t^{n+1-},\cdot))$ $= \sum_j|w_{j+1}^n-w_j^n|$ $= TV(w_{\Delta x}(t^n,\cdot))$. 
Denote $\calP_{\Delta x}$ the operator defined (for $x\in (x_{j-1/2},x_{j+1/2})$) by
$w_{\Delta x}(t^{n+1},x)$ $= \calP_{\Delta x}(w_{\Delta x}(t^{n+1-},x))$
$\equiv {1 \over \Delta x}\int_{x_{j-1/2}}^{x_{j+1/2}}w(t^{n+1-},y)\, dy$. 
This averaging operator $\calP_{\Delta x}$ is total variation diminishing, and so
$
 TV(w_{\Delta x}(t^{n+1},\cdot)) \leq TV(w_{\Delta x}(t^n,\cdot))$. $\hfill\Box$


\section{Numerical experiments}

\subsection{Formulation of the problem} 

In this section, coupling problems are numerically studied. Resonance in the infinitely sharp regime is closely investigated for smooth and discontinuous coupled solutions. The color function modeling the thick interface is chosen to be 
$v_\eta(x) = ( \erf(x/\eta)+1)/2$ ($x\in\RR$) 
for some fixed $\eta$. Here, $x\mapsto\erf(x)$ denotes the classical error function. Obviously $v_\eta$ takes value in the expected range $(0,1)$ and belongs to $W^{2,\infty}(\mathbb{R})$ (in fact to $\mathcal{C}^\infty(\mathbb{R})$). Observe that $x\mapsto v_\eta(x)-\frac{1}{2}$ is an even function so that the coupling formula (\ref{C0}) is, in some sense, ``symmetric'' 
in the left-- and right--hand partial differential equation.
In order to investigate the dependency of the discrete solutions upon the regularization, we consider the 
(two parameters) function 
$v_{\eta,\zeta}(x) = (\erf(x/\eta+\zeta)+1)/2$ ($x \in \RR$) 
with $\eta>0$ and $\zeta\in\RR$. Here, $\eta$ clearly monitors the thickness of the handshake coupling zone while $\zeta$ acts roughly speaking as a shifting perturbation  breaking the symmetry of the treatment of the left and right problems. All the calculations are performed over the bounded domain $[-1,1]$ and Neumann boundary conditions are used.

\subsection{Non--uniqueness for a resonant coupling problem}
\label{sec44} 

\begin{figure}[!ht]
\centering
\includegraphics[scale=1.4]{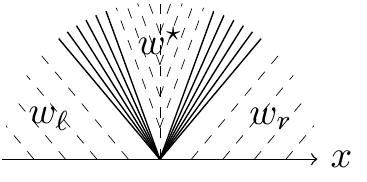}\\
$w_\star \in [w_\ell,w_r]$
\caption{A one-parameter family of smooth solutions in a resonant situation}
\label{fig-nonuniqueness}
\end{figure}

We 
consider the closure relations 
$f^-(w)=\frac{w^2}{2}$,
$f^+(w)=\frac{(w+1)^2}{2}$,
$\theta_-(w)=\theta_+(w)= w$, and 
the 
initial data  
is $w_0(x)= w_\ell=-1$ for $x<0$, and $w_r=3/2$ for $x>0$. 
The corresponding Riemann solution in the regime of an infinitely thin interface exhibits a resonance phenomenon which is depicted in Figure~\ref{fig-nonuniqueness}. Observe that in the present problem, 
the  sonic point 
associated to the right flux $f^+$ (respectively to the left flux $f^-$)
is $u^s_+=-1$ (resp. $u^s_-=0$). 
The ordering of the left-- and right--hand sonic states: 
$u^s_->u^s_+$
implies a resonance phenomenon in 
outgoing waves from the coupling interface. 
This is a sign of failure of uniqueness in solutions of the proposed coupled Riemann problem.  
We refer the reader to \cite{BoutinCoquelLeFloch09b}, the references therein, and to Figure~\ref{fig-nonuniqueness}.

\begin{figure}[!ht]
\centering
{\hfill
\subfloat[Numerical color functions]{\includegraphics[trim=60 60 30 60,clip,width=0.49\linewidth]{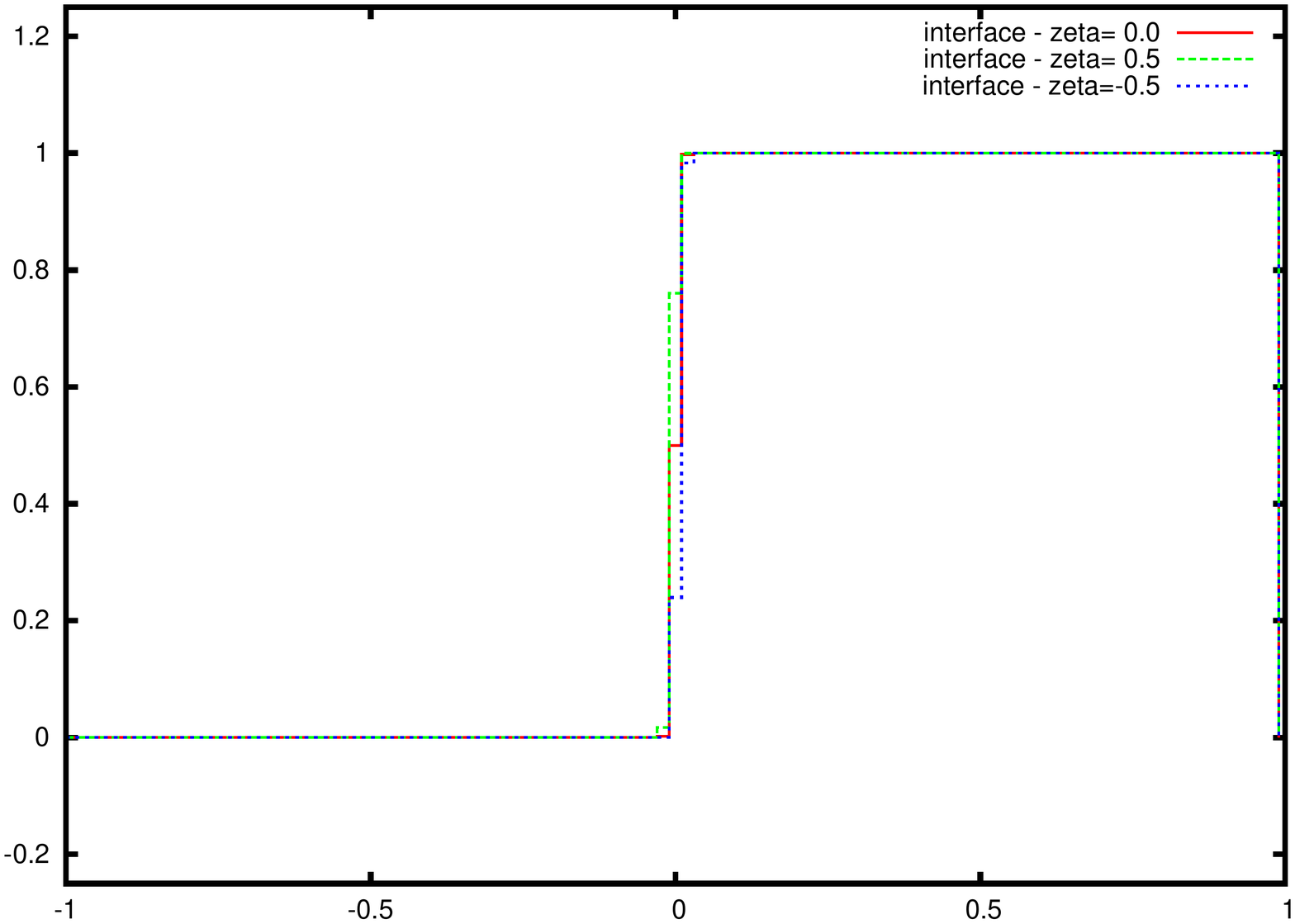}}
\hfill
\subfloat[Numerical solutions]{\includegraphics[trim=60 60 30 60,clip,width=0.49\linewidth]{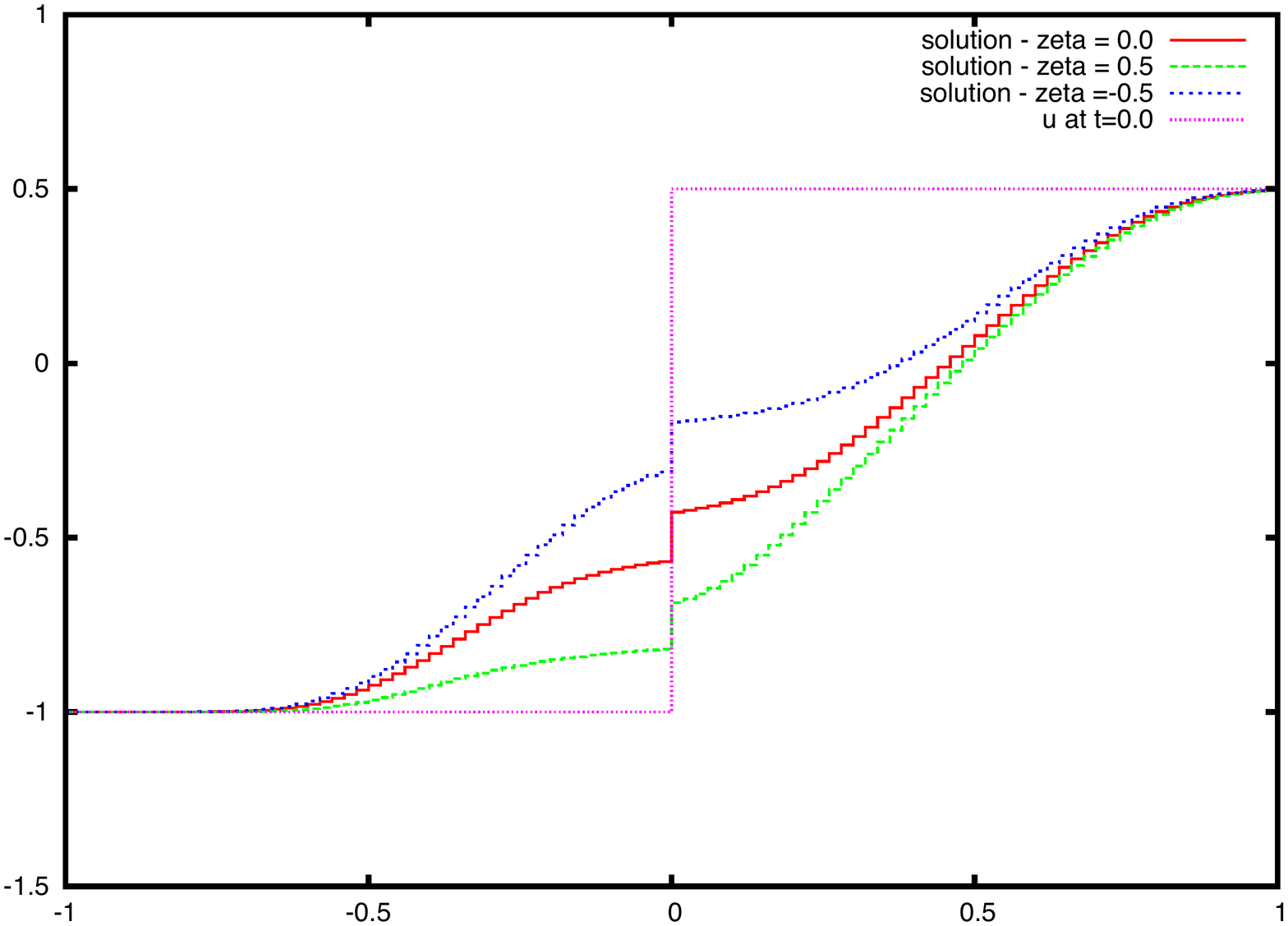}}
\hfill}
\caption{Failure of uniqueness in a resonant situation -- $N=100$}
\label{testDa}
\end{figure}

\begin{figure}[!ht]
\centering
{\hfill
\subfloat[Numerical color functions]{\includegraphics[trim=60 60 30 60,clip,width=0.49\linewidth]{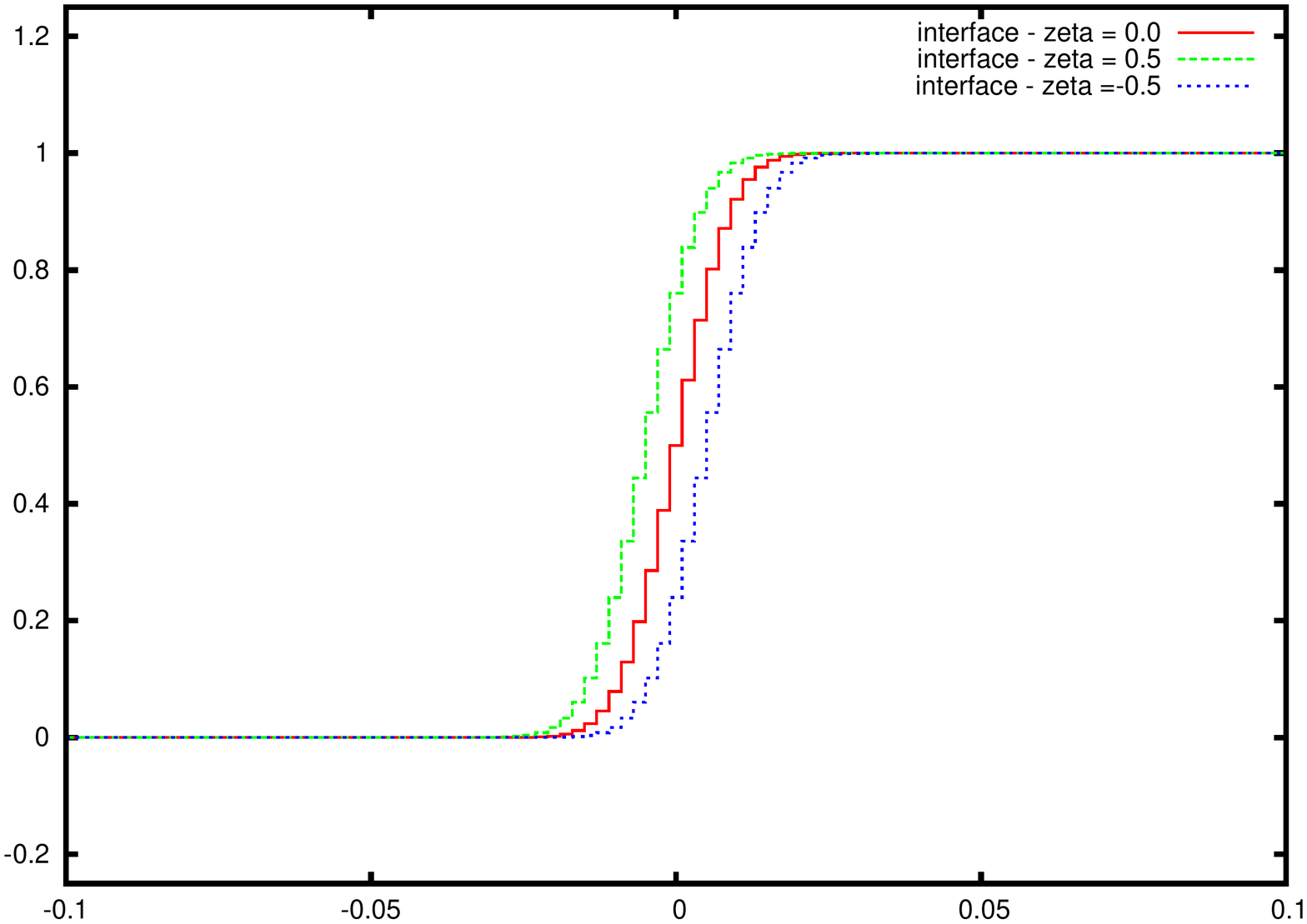}}
\hfill
\subfloat[Numerical solutions]{\includegraphics[trim=60 60 30 60,clip,width=0.49\linewidth]{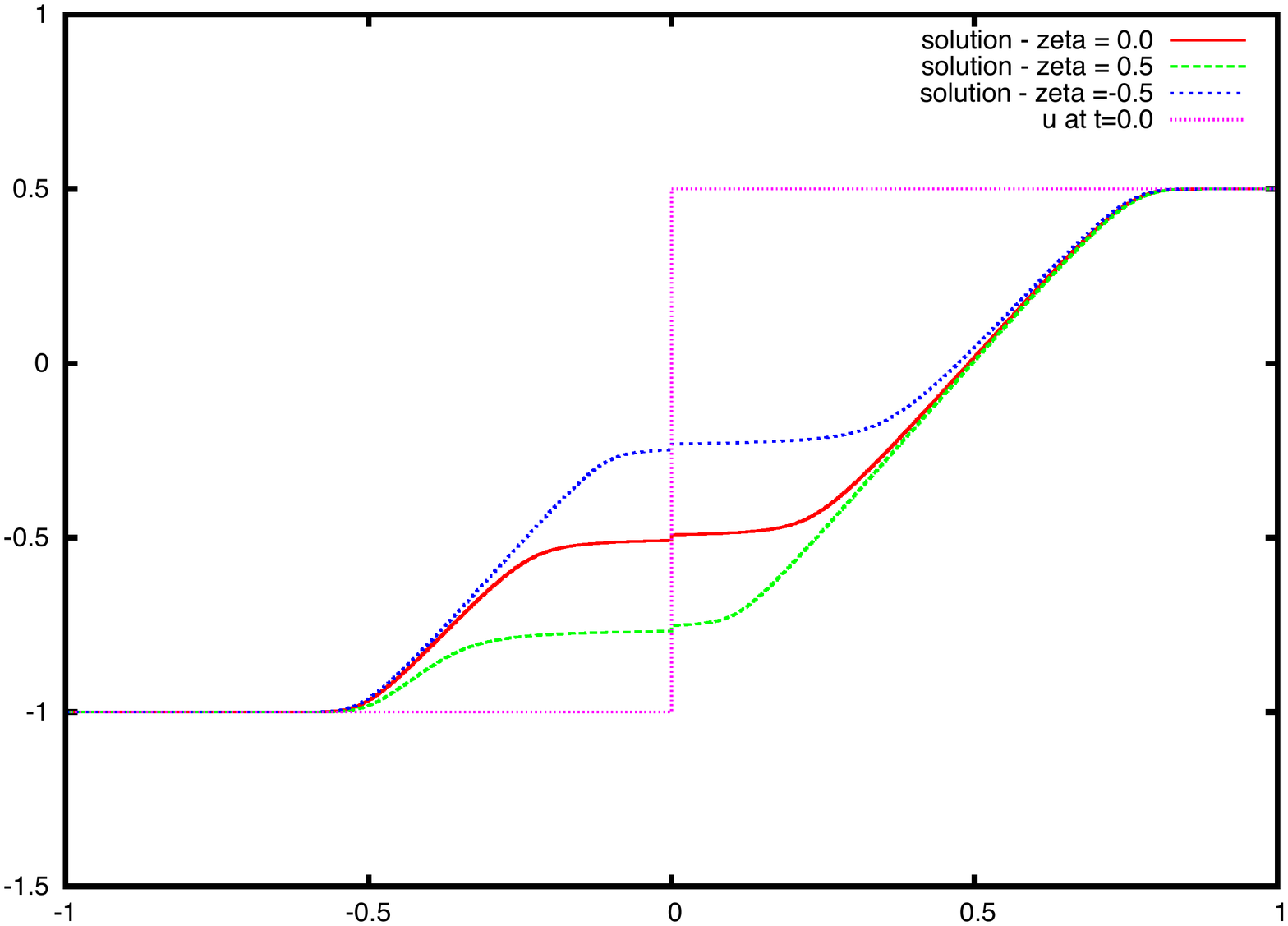}}
\hfill}
\caption{Failure of uniqueness in a resonant situation -- $N=1000$}
\label{testDb}
\end{figure}

\begin{figure}[!ht]
\centering
{\hfill
\subfloat[Numerical color functions]{\includegraphics[trim=60 60 30 60,clip,width=0.49\linewidth]{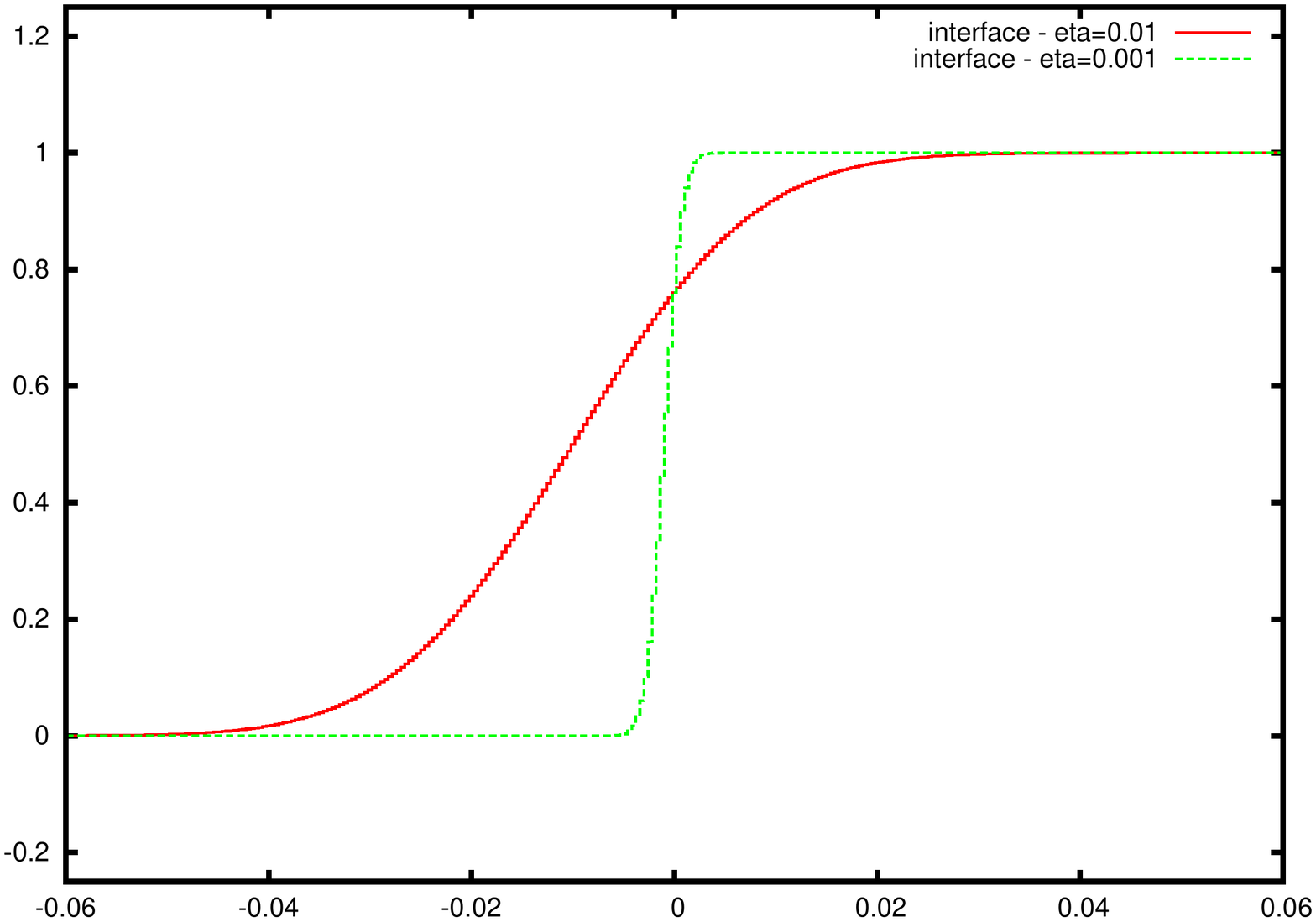}}
\hfill
\subfloat[Numerical solutions]{\includegraphics[trim=60 60 30 60,clip,width=0.5\linewidth]{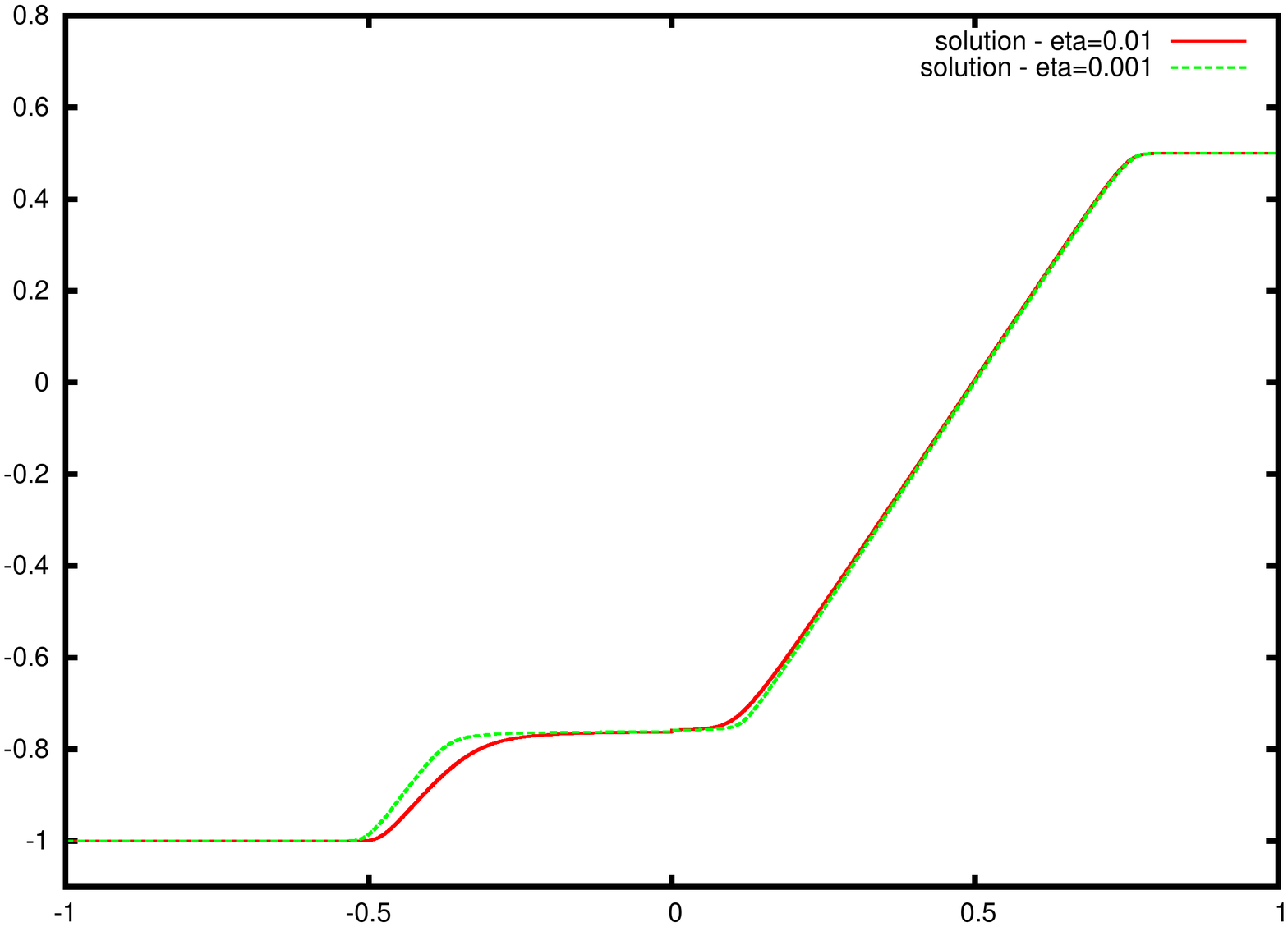}}
\hfill}
\caption{Sensitivity with respect to the thickness parameter}
\label{testE}
\end{figure}

In the regime of a regularized thick coupling interface, uniqueness of the solution of the Cauchy coupled problem is restored. A given regularization actually selects a given intermediate state so as to define a unique solution. Our aim here is to provide numerical evidences supporting this claim. We use the two-parameters family of regularizations 
$v_{\eta,\zeta}$ described above, for various thickness $\eta$ and shift $\zeta$.
Discrete solutions have been computed under the CFL condition (\ref{cfl}) using respectively 100 and 1000 grid points. The regularization parameters are respectively set to the constant value $\eta=5.10^{-3}$ and to three different values concerning $\zeta$: respectively $-0.5, 0$ and $0.5$. The extreme values select the left-- and right--hand problems depending of $\mbox{sgn}(\zeta)$ (negative values select $f^-$), while the intermediate value $0$ keeps a symmetry within the two problems.
The numerical results are displayed in Figures~\ref{testDa} and~\ref{testDb}. Two facts must be highlighted: the mesh refinement must be fine enough so as to capture a constant intermediate state, then and clearly the resulting value of the intermediate state $w_\star$ is highly sensitive with respect to the choice of the shifting parameter $\zeta$.  
In Figure~\ref{testE} we investigate the dependence with respect to the thickness parameter $\eta$, while the shifting parameter $\zeta$ is kept fixed at the constant value $0.5$. Two values of $\eta$ are considered: $\eta=0.01$ and $\eta=0.001$. The resulting $w$ profiles obtained using $5000$ grid points turn significantly less sensitive to the choice of the thickness parameter $\eta$ than to the choice of the shifting parameter $\zeta$.\\

\subsection{Another example of non--uniqueness}

\begin{figure}[!ht]
\centering
{\hfill
\subfloat[A left shock\label{fig-solshock-left}]{\includegraphics[scale=1]{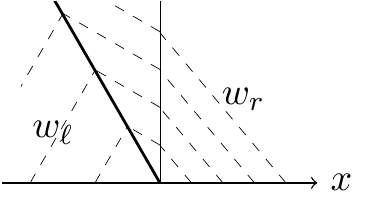}}%
\hfill
\subfloat[A standing discontinuity\label{fig-solshock-midd}]{\includegraphics[scale=1]{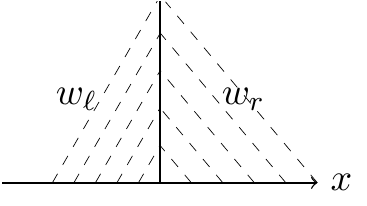}}%
\hfill
\subfloat[A right shock\label{fig-solshock-right}]{\includegraphics[scale=1]{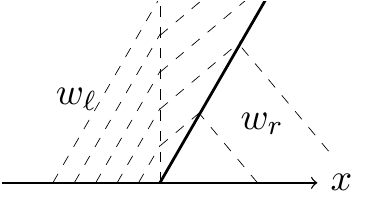}}%
\hfill}
\caption{Several discontinuous solutions in a resonant situation}
\label{fig-solutionconvergence2}
\end{figure}

\begin{figure}[!ht]
\centering
{\hfill
\subfloat[Numerical color functions]{\includegraphics[trim=60 60 30 60,clip,width=0.49\linewidth]{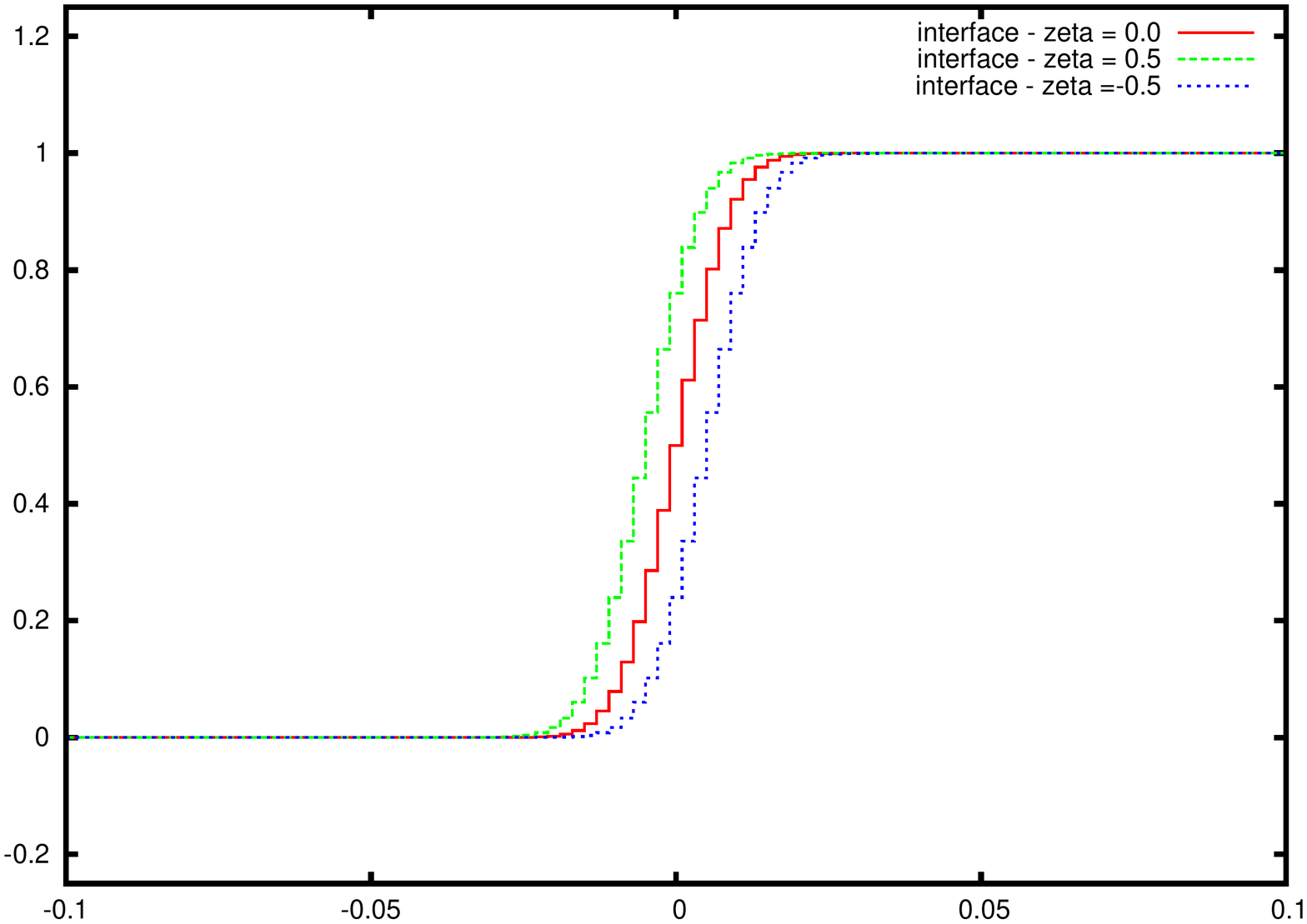}}
\hfill
\subfloat[Numerical solutions]{\includegraphics[trim=60 60 30 60,clip,width=0.49\linewidth]{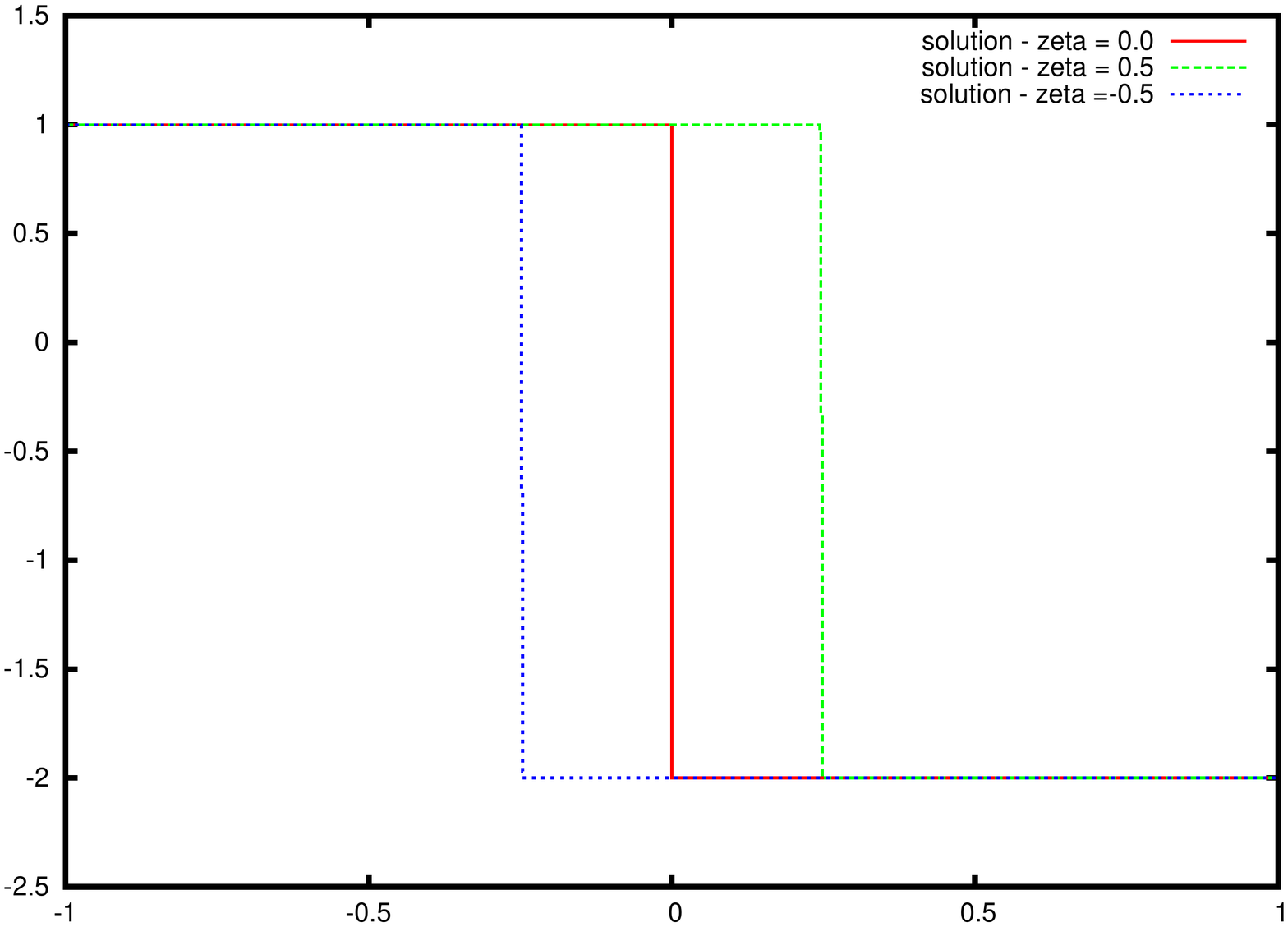}}
\hfill}
\caption{Failure of uniqueness in a resonant situation with discontinuous solutions -- $N=1000$}
\label{testF}
\end{figure}

This last test case is based on the closure relations in Section~\ref{sec44}, but
the initial data $w_0(x)=w_\ell=1$ for $x<0$ and $w_r=-2$ for $x>0$. 
Such a choice again results in a resonance phenomenon involving multiple solutions. Three different discontinuous solutions
are available in the regime of coupling interface with zero thickness. We refer the reader to \cite{BoutinCoquelLeFloch09a,BoutinCoquelLeFloch09b} for more details. These discontinuous solutions, depicted in Figure~\ref{fig-solutionconvergence2}, respectively coincide with a left moving shock, namely a discontinuity satisfying the Rankine-Hugoniot jump relation for $f^-$, a standing discontinuity, i.e.~a discontinuity with zero speed, and a right moving shock, that is a discontinuity satisfying the Rankine-Hugoniot jump relation for $f^+$. We show hereafter that these three different solutions are actually {\sl stable,}
 in the sense that we can capture each of them numerically 
by choosing suitable regularization parameters in the thick coupling interface $v_{\eta,\zeta}$. We choose as previously a fixed thickness $\eta=5.10^{-3}$ and we make the shifting parameter $\zeta$ to vary from the negative value $-0.5$ (selecting the left--hand equation) to the positive one $0.5$ (selecting the right--hand equation) and the vanishing value $\zeta=0$, keeping the symmetry in between the two. As heuristically expected, the results displayed in Figure~\ref{testF} show that the negative value of $\zeta$ captures the left moving shock, the right--hand state selects the right moving shock, while $0$ restores the standing discontinuity. In particular, multiple solutions may be stable.

\vskip-.5cm

\begin{thebibliography}{10}

\bibitem{AdimurthiMishraGowda05}
{\sc Adimurthi, S. Mishra, and G.D.V. Gowda,}
Optimal entropy solutions for conservation laws with discontinuous flux-functions, 
J. Hyperbolic Differ. Equ. 2 (2005), 783--837.

\bibitem{Amadori08}
{\sc D. Amadori, L. Gosse, G. Graziano,}
Godunov-type approximation for a general resonant balance law with large data,
 J. Differential Equations 198 (2004), 233--274. 
 
\bibitem{AmbrosoEtAl07}
{\sc A. Ambroso, C. Chalons, F. Coquel, E. Godlewski, F. Lagouti\`ere, P.-A. Raviart, and N. Seguin,}
A relaxation method for the coupling of systems of conservation laws,
in ``Hyperbolic problems: Theory, Numerics, Appl.'', Springer Verlag, Berlin, 2008, pp.~947--954.

\bibitem{AmbrosoEtAl08}
{\sc A. Ambroso, C. Chalons, F. Coquel, E. Godlewski, F. Lagouti\`ere, P.-A. Raviart, and N. Seguin,}
Coupling of general Lagrangian systems, 
Math. Comp. 77 (2008), 909--941.

\bibitem{AmbrosoEtAl08b}
{\sc A. Ambroso, C. Chalons, F. Coquel, E. Godlewski, F. Lagouti\`ere, P.-A. Raviart, and N. Seguin,}
Relaxation methods and coupling procedures, 
Internat. J. Numer. Methods Fluids 56 (2008), 1123--1129.

\bibitem{AmbrosoEtAl08c}
{\sc A. Ambroso, C. Chalons, F. Coquel, T. Gali\'e, E. Godlewski, P.-A. Raviart, and N. Seguin,} 
The drift-flux asymptotic limit of barotropic two-phase two-pressure models,
Commun. Math. Sci. 6 (2008), 521--529.

\bibitem{AmbrosoEtAl07b}
{\sc A. Ambroso, C. Chalons, F. Coquel, E. Godlewski, F. Lagouti\`ere, P.-A. Raviart, and N. Seguin,}
The coupling of homogeneous models for two-phase flows, 
Int. J. Finite Vol. 4 (2007), 39--54.

\bibitem{ALO} {\sc P. Amorim, P.G. LeFloch, and B. Okutmustur,} 
Finite volume schemes on Lorentzian manifolds, 
Comm. Math. Sc. 6 (2008), 1059--1086. 

\bibitem{AudussePerthame05}
{\sc E. Audusse and B. Perthame,}
Uniqueness for scalar conservation laws with discontinuous flux via adapted entropies,
Proc. Roy. Soc. Edinburgh Sect. A 135 (2005), 253--265.

\bibitem{BachmannVovelle06}
{\sc F. Bachmann and J. Vovelle,} 
Existence and uniqueness of entropy solution of scalar conservation laws with a flux function involving discontinuous coefficients,
Comm. Partial Differential Equations 31 (2006), 371--395.

\bibitem{neptune}
{\sc D. Bestion, M. Boucker, P. Boudier, P. Fillion, M. Grandotto, A. Guelfi, J.M. H{\'e}rard, E. Hervieu, and P. P{\'e}turaud,} 
Neptune: a new software platform for advanced nuclear thermal hydraulics, Nuclear Sc. Eng. 156 (2007), 281--324.

\bibitem{BenArtziLeFloch07}
{\sc M. Ben-Artzi and P.G. LeFloch,}
Well-posedness theory for geometry-compatible hyperbolic conservation laws on manifolds, 
Ann. Inst. H. Poincar\'e -- Anal. Nonlin\'eaire 24 (2007), 989--1008.

\bibitem{Bouchut04} 
{\sc F. Bouchut,} 
{\em Nonlinear stability of finite volume methods for hyperbolic conservation laws and well-balanced schemes for sources,}
``Frontiers in Mathematics'', Birkh\"auser Verlag, B\"asel, 2004.

\bibitem{BoutinChalonsRaviart10}
{\sc B. Boutin, C.~Chalons, and P.-A. Raviart,} 
Existence result for the coupling problem of two scalar conservation laws with Riemann initial data,
Math. Models Methods Appl. Sci. 20 (2010), 1859--1898.

\bibitem{BoutinCoquelGodlewski08}
{\sc B. Boutin, F. Coquel, and E. Godlewski,} 
Dafermos' regularization for interface coupling of conservation laws,
in ``Hyperbolic problems: Theory, Numerics, Applications'', Springer Verlag, Berlin, 2008, pp.~567--575.

\bibitem{BoutinCoquelLeFloch09a}
{\sc B. Boutin, F. Coquel, and P.G. LeFloch,} 
Coupling techniques for nonlinear hyperbolic equations. I. Self-similar diffusion for thin interfaces, 
Proc. Roy. Soc. Edinburgh Sect. A 141 (2011), 921--956.

\bibitem{BoutinCoquelLeFloch09b}
{\sc B. Boutin, F. Coquel, and P.G. LeFloch,} 
Coupling techniques for nonlinear hyperbolic equations. II, in preparation. 
 
\bibitem{BoutinCoquelLeFloch09d}
{\sc B. Boutin, F. Coquel, and P.G. LeFloch,} 
Coupling techniques for nonlinear hyperbolic equations. 
IV.  Multi--component coupling and multidimensional well--balanced schemes, 
Preprint ArXiv:1206.0248. 

\bibitem{BurgerKarlsen08}
{\sc R. B\"urger and K.H. Karlsen,} 
Conservation laws with discontinuous flux: a short introduction,
J. Engrg. Math. 60 (2008), 241--247.

\bibitem{BurgerKarlsenTowers09}
{\sc R. B\"urger, K. H. Karlsen, and J. D. Towers,}
{An Engquist-Osher-type scheme for conservation laws with discontinuous flux adapted to flux connections,}
SIAM J. Numer. Anal. 47 (2009) 1684--1712.

\bibitem{ChainaisCHampier01}
{\sc C. Chainais--Hillairet and S. Champier,} 
Finite volume schemes for nonhomogeneous scalar conservation laws: error estimate,
Numer. Math. 88 (2001), 607--639.

\bibitem{ChalonsRaviartSeguin08}
{\sc C. Chalons, P.-A. Raviart, and N. Seguin,} 
The interface coupling of the gas dynamics equations,
Quart. Appl. Math. 66 (2008), 659--705.

\bibitem{Coquel12}
{\sc F. Coquel,}
Coupling of nonlinear hyperbolic systems :
A journey from mathematical to numerical issues, in ``Numerical methods for hyperbolic equations'', 
E.  V{\'a}zquez-Cend{\'o}n, A. Hidalgo, P. Garcia-Navarro, and L. Cea editors, Taylor and Francis Group, 2012.

\bibitem{CoquelLeFloch93}
{\sc F. Coquel and P.G. LeFloch,}
Convergence of finite difference schemes for conservation laws in several space dimensions: a general theory,
SIAM J. Numer. Anal. 30  (1993), 675--700.

\bibitem{Dafermos73}
{\sc C.M. Dafermos,} 
Solution of the {R}iemann problem for a class of hyperbolic systems of conservation laws by the viscosity method,
Arch. Rational Mech. Anal. 52 (1973), 1--9.

\bibitem{DLM}
{\sc G. Dal Maso, P.G. LeFloch, and F. Murat,}
Definition and weak stability of nonconservative products, 
J. Math. Pures Appl. 74 (1995), 483--548. 

\bibitem{DiPerna85}
{\sc R.J. DiPerna,} 
Measure--valued solutions to conservation laws,
Arch. Rational Mech. Anal. 88 (1985), 223--270.

\bibitem{DuboisLeFloch88}
{\sc F. Dubois and P.G. LeFloch,} 
Boundary conditions for nonlinear hyperbolic systems of conservation laws, 
 J. Differential Equations 71 (1988), 93--122.

\bibitem{Duret04}
{\sc E. Duret, Y. Peysson, Q. H. Tran, and P. Rouchon,}
Active bypass to eliminate severe-slugging in multiphase production, 
Proceedings of the 4th North American Conference on Multiphase Technology, Ban, June 2004, J. Brill, ed.,
Cranfield, 2004, BHR Group Limited, 205--223.

\bibitem{EymardGallouetHerbin00}
{\sc R. Eymard, T. Gallou\"et, and R. Herbin,} 
Finite volume methods, in ``Handbook of numerical analysis'',
 Vol.~VII, North--Holland, Amsterdam, 2000, 713--1020.

\bibitem{GadelHak}
{\sc M. Gad-El-Hak,}
{\em Flow control: Passive, active, and reactive flow management,}
Cambridge University Press, 2006.

\bibitem{GoatinLeFloch04}
{\sc P. Goatin and P.G. LeFloch,} 
The Riemann problem for a class of resonant hyperbolic systems of balance laws, 
Ann. Inst. H. Poincar\'e Anal. Non Lin\'eaire 21 (2004), 881--902. 

\bibitem{GodlewskiRaviart04}
{\sc E. Godlewski and P.-A. Raviart,} 
The numerical interface coupling of nonlinear hyperbolic systems of conservation laws. I. The scalar case,
Numer. Math. 97 (2004), 81--130.

\bibitem{Gosse01}
{\sc L. Gosse,} 
A well-balanced scheme using non-conservative products designed for hyperbolic systems of conservation laws with source terms, 
Math. Models Methods Appl. Sci. 11 (2001), 339--365. 

\bibitem{GreenbergLeroux92}
{\sc J.M. Greenberg and A.Y. Leroux,} 
A well-balanced scheme for the numerical processing of source terms in hyperbolic equations,
SIAM J. Numer. Anal. 33 (1996), 1--16. 

\bibitem{HL} 
{\sc T.Y. Hou and P.G. LeFloch,}  
Why nonconservative schemes converge to wrong solutions. Error analysis, 
Math. of Comput. 62 (1994), 497--530.

\bibitem{IsaacsonTemple92}
{\sc E. Isaacson and B.J. Temple,}
Nonlinear resonance in systems of conservation laws,
SIAM J. Appl. Math. 52 (1992), 1260--1278.

\bibitem{JL}
{\sc K.T. Joseph and P.G. LeFloch,}
Boundary layers in weak solutions of hyperbolic conservation laws. II. Self-similar vanishing diffusion limits, 
Comm. Pure Appl. Anal. 1 (2002), 51--76.

\bibitem{Kroener98}
{\sc D.~Kr\"oner,} 
Finite volume schemes in multidimensions,
in ``Numerical analysis 1997'' (Dundee), Pitman Res. Notes Math. Ser.,
 Vol. 380, Longman, Harlow, 1998, pp.~179--192.

\bibitem{Kruzkov70}
{\sc S.N. Kruzkov,} 
First order quasilinear equations with several independent variables, 
Mat. Sb. (N.S.) 81 (1970), 228--255.

\bibitem{LeFloch88}
{\sc P.G. LeFloch,} 
Entropy weak solutions to nonlinear hyperbolic systems in nonconservative form, 
Comm. Part. Diff. Equa. 13 (1988), 669--727.

\bibitem{LeFloch89}
{\sc P.G. LeFloch,} 
Shock waves for nonlinear hyperbolic systems in nonconservative form,
Institute for Math. and its Appl., Minneapolis, Preprint \# 593 (1989). Available at
www.ima.umn.edu/preprints/Jan89Dec89/593.pdf.


\bibitem{LeFloch-interface}
{\sc P.G. LeFloch,}
An existence and uniqueness result for two nonstrictly hyperbolic systems, 
IMA Vol. in Math. and its Appl.,``Nonlinear evolution equations that change type'', 
ed.~B.L. Keyfitz and M. Shearer, Springer Verlag, Vol. 27, 1990, 126--138.

\bibitem{LeFloch-ARMA}
{\sc P.G. LeFloch,}
Propagating phase boundaries. Formulation of the problem and existence via the Glimm scheme, 
Arch. Rational Mech. Anal. 123 (1993), 153--197.

\bibitem{LeFloch-book}
{\sc P.G. LeFloch,}
{\it Hyperbolic systems of conservation laws. The theory of classical and nonclassical shock waves,} 
Lectures in Mathematics, ETH Z\"urich, Birkh\"auser, 2002. 

\bibitem{LT}
{\sc P.G. LeFloch and M.-D. Thanh,}  
A Godunov-type method for the shallow water equations with variable topography in
 the resonant regime,  J. Comput. Phys. 230 (2011), 7631--7660.

\bibitem{LinSchecter03}
{\sc X.-B. Lin and S. Schecter}, 
Stability of self-similar solutions of the Dafermos regularization of a system of conservation laws,
SIAM J. Math. Anal. 35 (2003), 884--921.

\bibitem{SchecterPlohrMarchesin04}
{\sc S. Schecter, B.J. Plohr, and D. Marchesin,}  
Computation of Riemann solutions using the Dafermos regularization and continuation,
Discrete Contin. Dyn. Syst. 10 (2004), 965--986.

\bibitem{SeguinVovelle03}
{\sc N. Seguin and J. Vovelle,} 
Analysis and approximation of a scalar conservation law with a flux function with discontinuous coefficients,
Math. Models Methods Appl. Sci. 13 (2003), 221--257.

\bibitem{Szepessy89}
{\sc A. Szepessy,} 
An existence result for scalar conservation laws using measure valued solutions,
Comm. Partial Differential Equations 14 (1989), 1329--1350. 

\end{thebibliography}
\end{document}